\documentclass[%
 aip,
 amsmath,amssymb,
 reprint,%
]{revtex4-1}

\usepackage{graphicx}
\usepackage{dcolumn}
\usepackage{bm}

\usepackage[utf8]{inputenc}
\usepackage[T1]{fontenc}
\usepackage{mathptmx}
\usepackage{etoolbox}

\usepackage[english]{babel}
\usepackage{amsmath}
\usepackage{systeme} 
\usepackage{amsfonts} 

\newcommand{\ignore}[2]{\hspace{0in}#2}
\usepackage[colorlinks=true, allcolors=blue]{hyperref}

\usepackage{caption}

\usepackage{floatrow}
\usepackage{subfig}
\newsavebox{\measurebox}

\usepackage{multirow}

\makeatletter
\def\@email#1#2{%
 \endgroup
 \patchcmd{\titleblock@produce}
  {\frontmatter@RRAPformat}
  {\frontmatter@RRAPformat{\produce@RRAP{*#1\href{mailto:#2}{#2}}}\frontmatter@RRAPformat}
  {}{}
}%
\makeatother
\begin{document}

\preprint{AIP/123-QED}

\title[Ordinal Poincar\'e Sections: Reconstructing the First Return Map from an Ordinal Segmentation of Time Series]{Ordinal Poincar\'e Sections: Reconstructing the First Return Map from an Ordinal Segmentation of Time Series}
\author{Zahra Shahriari}
\author{Shannon Dee Algar}%
\author{David M. Walker}%
\author{Michael Small}
 \altaffiliation[Also at ]{Mineral Resources, CSIRO, Kensington, Western Australia 6151, Australia}
 \email{zahra.shahriari@research.uwa.edu.au}
\affiliation{ 
Complex Systems Group, Department of Mathematics and Statistics, The University of Western Australia, Crawley, Western Australia 6009, Australia
}%


\date{\today}

\begin{abstract}
We propose a robust and computationally efficient algorithm to generically construct first return maps of dynamical systems from time series without the need for embedding. Typically, a first return map is constructed using a heuristic convenience (maxima or zero-crossings of the time series, for example) or a computationally delicate geometric approach (explicitly constructing a Poincar\'e section from a hyper-surface normal to the flow and then interpolating to determine intersections with trajectories). Our approach relies on ordinal partitions of the time series and builds the first return map from successive intersections with particular ordinal sequences. Generically, we can obtain distinct first return maps for each ordinal sequence. We define entropy-based measures to guide our selection of the ordinal sequence for a ``good'' first return map and show that this method can robustly be applied to time series from classical chaotic systems to extract the underlying first return map dynamics. The results are shown on several well-known dynamical systems (Lorenz, R{\"o}ssler and Mackey-Glass in chaotic regimes).

\end{abstract}

\maketitle

\begin{quotation}
Nonlinear time series analysis seeks to uncover features of a deterministic dynamic system from observed (usually scalar) time series data. Methods have been developed to build models of the underlying vector field, estimate dynamical invariants, or build a chart of the geometry of the attractor and its unstable periodic orbits. The first return map, which is one of the most fundamental descriptions of the dynamics of a dynamical system, is typically estimated by developing a complex numerical model of the flow or by employing heuristic approaches and trying to find the best. Nonetheless, first return maps remain one of the simplest and most straightforward mechanisms with which to understand complicated dynamics --- whenever they can be computed. In particular, for three-dimensional chaotic systems, the first return map is a map from $\mathbb{R}$ to $\mathbb{R}$, and the interesting behavior is often encapsulated in the behavior of familiar "one-hump" maps.

We propose a generic alternative based on the recent increased interest in the ordinal encoding of dynamics and symbolic dynamics in general. Our method is useful as it provides a straightforward, robust, generic, and numerically simple approach to access the first return map of a continuous dynamical system from a scalar time series. We show that this method can successfully reconstruct the first return map for simple low-dimensional chaotic systems as well as infinite-dimensional delay differential systems. This method provides a simple approach to directly identify the evolution of chaotic dynamics in experimental systems from measured time series. 


\end{quotation}

\section{\label{sec:level1}Introduction}

The Poincar\'e map (also known as the first return map) is a useful tool for studying dynamical systems --- projecting a continuous dynamical system to a lower dimensional map. A Poincar\'e section is the intersection of an orbit in the state space of a continuous dynamical system with a particular lower-dimensional surface transverse to the system's flow \cite{r33}. More specifically, consider an orbit with initial conditions within a section of space, which then leaves that section, and observe the point at which this orbit first returns to the section.
The first return map (FRM) is then created by making a map that connects the first point in the section to the second, the second point to the third, and so on. The transversality of the Poincar\'e section means that the orbits begin in the subspace and flow through it rather than parallel to it. Hereafter we will refer to the $n-1$ dimensional surface transverse to the flow as the {\em Poincar\'e section}, and the map we construct from successful intersections with that section as the {\em First Return Map}.

Kakutani introduced the FRM in 1943 following the idea of the Poincar\'e recurrence theorem with the aim of investigating the properties of induced measure-preserving transformations \cite{r31}. A Poincar\'e map can be considered as a discrete dynamical system with dimension $n-1$ from the original continuous dynamical system of dimension $n$. It is frequently used to analyze the dynamical system more simply because it preserves many properties of the original system transforming periodic and quasiperiodic orbits into periodic and quasiperiodic points. The FRM has a lower-dimensional state space and is a discrete map. 

\begin{figure}[htbp]
    \centering
    \sidesubfloat[]
        {\includegraphics[width = 0.9\linewidth]{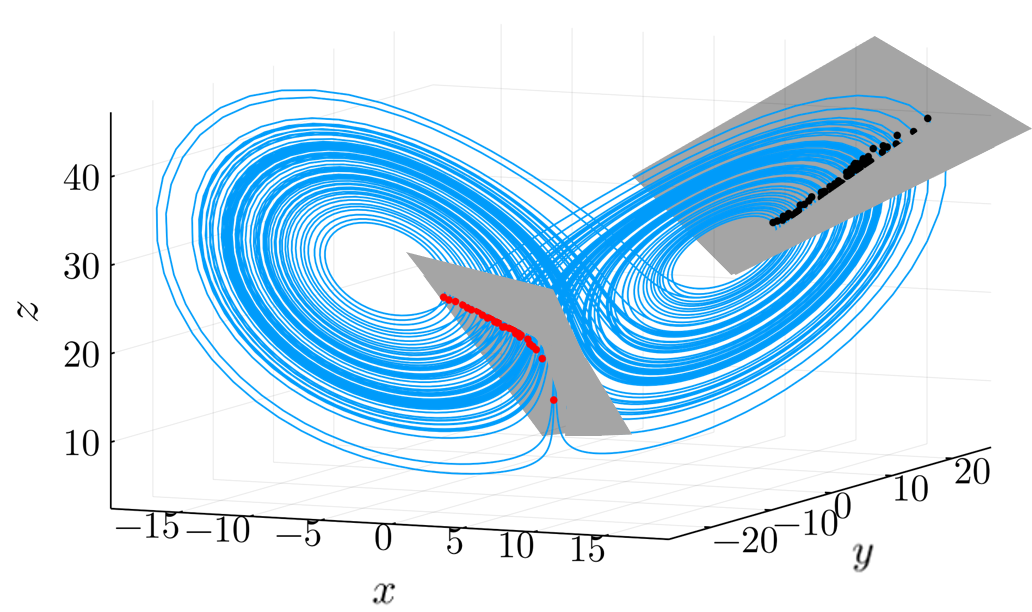}
        \label{LorenzPoincareSection:a}}\\[-1pt]
   \sidesubfloat[]
        {\includegraphics[width = 0.9\linewidth]{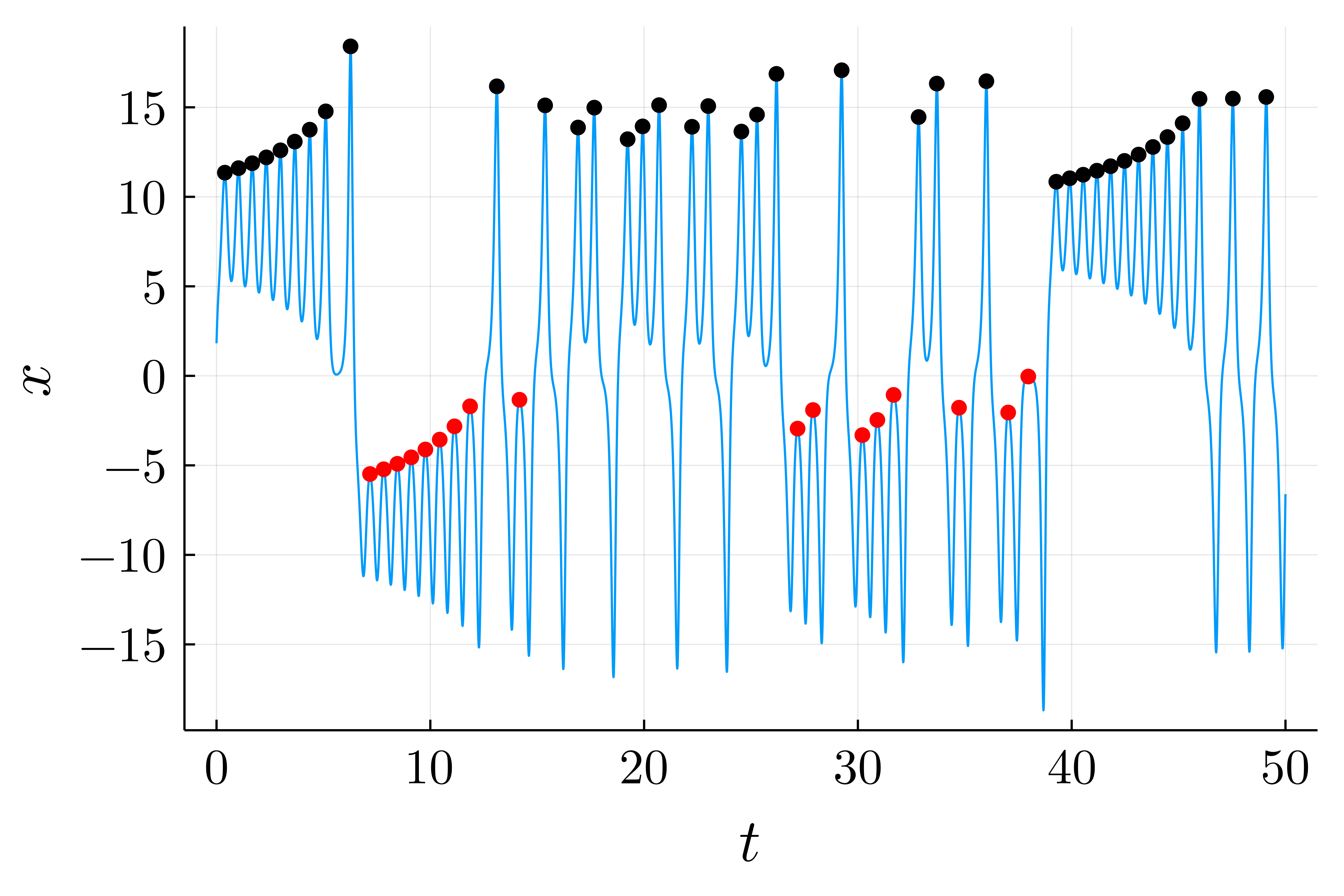}
        \label{LorenzPoincareSection:b}}\\[-1pt]
    \sidesubfloat[]
        {\includegraphics[width = 0.9\linewidth]{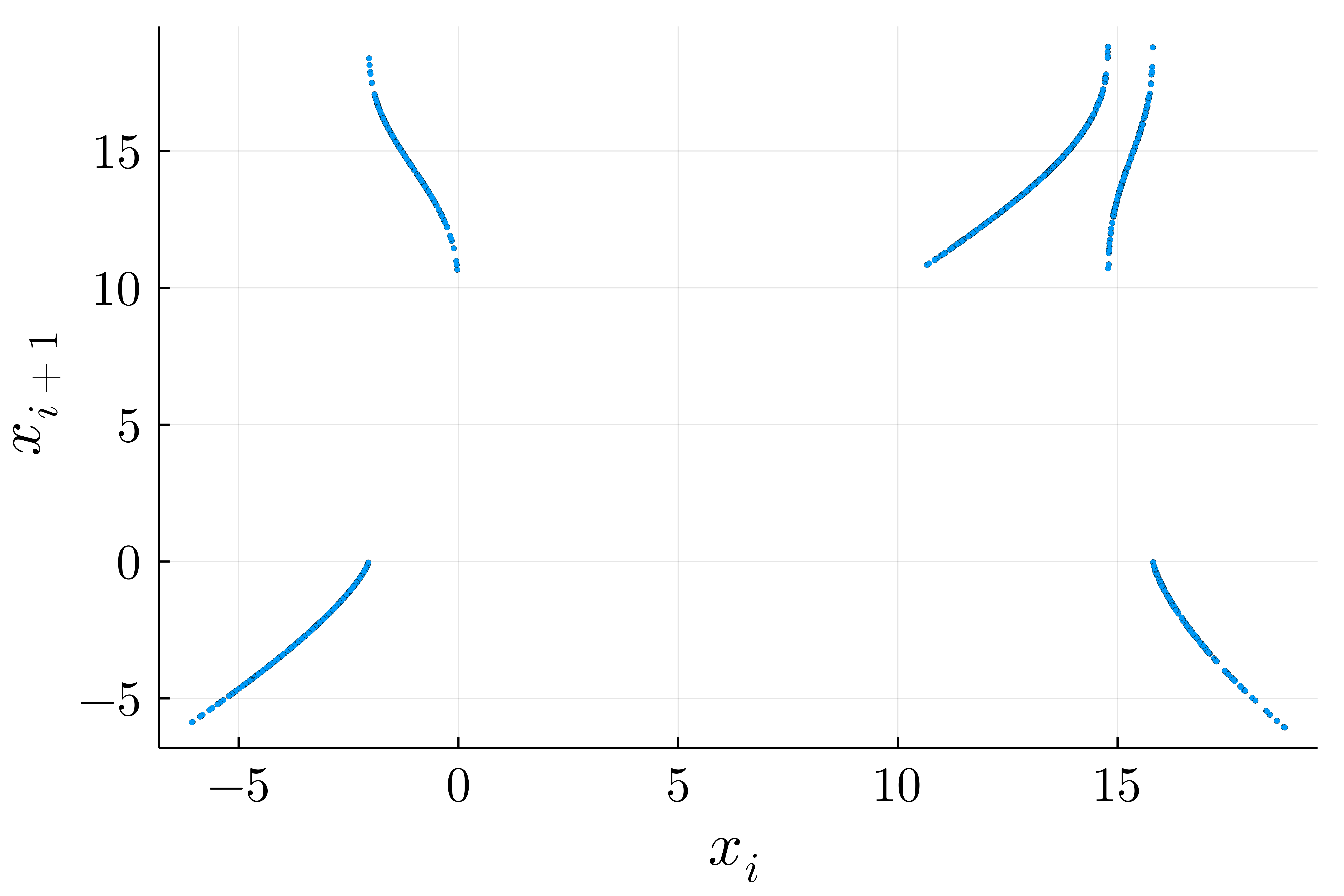}
        \label{LorenzPoincareSection:c}}
    \caption{(a) The surfaces of maxima of $x$ on the Lorenz attractor. The entry points to these surfaces are indicated by red and black dots. (b) Part of the Lorenz $x$ time series after removing the transient. The local maxima in the $x < 0$ region are depicted by red dots, whereas the local maxima in the $x > 0$ region are depicted by black dots. These two sets of points correspond to the points shown in Fig. 1a on the time series. (c) Corresponding values of relative maxima of $x$ and subsequent relative maxima of $x$, known as first return map (FRM).}
    
    \label{LorenzPoincareSection}
\end{figure}

These FRMs are valuable tools for studying oscillatory dynamics and have been used to explain the origins of the chaotic dynamics in specific models \cite{r18}. 
Figure~\ref{LorenzPoincareSection:a} shows the surfaces of maxima of $x$ with $\dot{x}=0$ in gray as a Poincar\'e section on the Lorenz attractor. The entering points to this section are marked with dots on the attractor.
It is clear that one of the $\dot{x}=0$ surfaces crosses the attractor in a region where $x$ is negative (red dots), and the other one passes the attractor where $x$ is positive (black dots).
Figure~\ref{LorenzPoincareSection:b} depicts a portion of the Lorenz time series after passing through the transient. Red and black dots highlight the local maxima entering the $\dot{x}=0$ Poincar\'e section. 
The FRM generated from the red and black dots is represented in Fig.~\ref{LorenzPoincareSection:c}. The idea of constructing the FRM from the points entering a Poincar\'e section was first introduced in a similar figure in 1963 for the Lorenz attractor's $z$ time series in Lorenz's famous paper \cite{r19}.

Despite all the desirable features of the first return map, finding an ideal Poincar\'e section that can correctly capture the system's characteristics is not always easy \cite{r44}. Also, suppose the system's equations are not fully available, and only its output time series is known. In that case, to use the Poincar\'e section, it is necessary to transform the time series into a suitable higher-dimensional structure that allows the application of standard methods for analyzing high dimension \cite{r22}.
One of the most popular methods to do so is the reconstruction of phase space from time series in a way that is equivalent to a system's original phase space \cite{r35}. The paradigmatic example of embedding techniques is the Takens' delay embedding theorem for an infinite, noise-free time series to unfold the attractor of its generator system \cite{r36}. This method is widely applied in practice; however, it is not always easy to find the appropriate dimension and lag in the embedding solution that produces an attractor topologically similar to the original system and is computationally efficient \cite{r43}.

In addition to embedding, recent approaches to convert time series data to a complex network form have gained popularity and allow one to examine time series in higher dimensions and use complex network analysis tools to extract time series features \cite{r29,r30,r41}. Several algorithms can generate a complex network representation from a time series \cite{r12,r13,r14,r15}. The one that we are going to use in this study is the ordinal partition network approach \cite{r16}. Converting a continuous time series (albeit sampled at discrete times) to a complex network involves two important approximations. Nodes of the networks are associated with regions of state space, and so there is a coarse-grained quantization of the original observations. Moreover, the network representation encodes transitions between nodes based on probability. It allows one to compute properties related to the underlying invariant measure but does not preserve the deterministic dynamics \cite{r23}.

One of the most common methods for generating a complex network from a time series is to form an ordinal partition network (OPN), which was introduced for the first time in 2013 \cite{r16}. 
This dynamical symbolic mapping \cite{r23}, which employs Markov modeling of dynamical systems, enables the approximation of deterministic dynamics by examining the statistical properties of a stochastic model. 
In simpler terms, this partitioning defines a map between the time delay embedded in observed data and the symbolic dynamics of interest. 

Unlike OPN approaches, which divide state space into distinct regions and provide a stochastic approximation to the dynamics, the ordinal-based FRM approach we propose here preserves the original fidelity while constructing the deterministic discrete map.
We construct it by assigning an ordinal partition to each point in the time series and taking the change in the ordinal partition between two nearby points as the entrance to that specific ordinal partition section.

In our paper, we present a generic method for detecting the dynamical behavior of a time series by treating the ordinal partitions as Poincar\'e sections. Although these sections differ slightly from the primary definition of Poincar\'e sections, we show that they can behave similarly and may be used to extract the main dynamical features of the system. This method extracts features from a single time series without requiring access to the entire attractor or performing a full embedding.
Furthermore, by combining the FRMs from different ordinal partitions, additional information regarding the dynamical features can be elucidated. It is also noise resistant as it is built on ordinal partitions, which allows it to be used on real-world data with significant noise contamination \cite{r12}.
We used a numerically generated time series from the Lorenz system in Eq.~\ref{Lorenz} to validate our theoretical discussion.  We assign a symbol to each point of the time series by ordinal partitioning the amplitude of the time series points in each window. For the next step, for each set of ordinal partitions, we create a multiscale ordinal partition network from all the points in it to calculate the entropy for each set. We can determine which ordinal partitions are good Poincar\'e sections by ranking them based on their OPN's entropies and comparing their FRMs to those obtained by the original forms of Poincar\'e sections.

The rest of this paper is organized as follows. We first introduce the methods and measures that have been used to encode the time series in a symbolic sequence that allows one to estimate the dynamical features of the time series in Sec. 2, where we also introduce different types of permutation entropies for ranking ordinal partitions. Section 3 presents the results of applying this method to the Lorenz time series, and Section 4 concludes. The results of this method's investigation on Lorenz time series with other ordinal parameters, as well as time series from some other dynamical systems, are given in the appendix.

\section{\label{sec:Methods}Methods}

\subsection{\label{sec:Ordinal Partition Networks}Ordinal Partition Networks}

Let $x_t$ represent a scalar time series measured from a known or unknown dynamical system. We generate an OPN from this time series in the following manner. We consider a window of length $L$ and slide it over the time series to generate an ordinal partition network (OPN). Each window is made up of $m$ points, with $\tau$-points gap between each pair. This means that the window size is $L = (m-1)\tau$. As a result, the corresponding ordinal partition $o^{(i)} = (\pi_1, \pi_2,..., \pi_m)$ where $\pi_j \in \{1, 2,...,m\}$ ($\pi_j \neq \pi_k \text{ if } j \neq k$) could be defined for each vector $z^{(i)} = (x_i, x_{i+\tau},..., x_{i+(m-1)\tau})$. 
In this way, $o^{(i)}$ captures the relative order of $z^{(i)}$ vector elements. 
This symbolic ordering can be done in two different ways: \textit{Amplitude Ranking}, and \textit{Chronological Index Ranking} \cite{r23}. Amplitude ranking ranks the points based on relative amplitude, so the point with the lowest amplitude has rank $m$ and the point with the highest amplitude has rank $1$. 
On the other hand, chronological index ranking uses the time index to sort the points such that it is the index of the lowest value that is placed in first in $o$. 
Table \ref{tab:Ranking} shows these two ways of ranking for $m=3$. In this work, we used chronological index ranking for creating the OPN.

\begin{table}    
\caption{\label{tab:Ranking}Amplitude ranking and chronological index ranking
for the window size $m=3$.}
\begin{ruledtabular}
\begin{tabular}{ccc}
Order of Points&Amplitude Ranking&Chronological\\
&&index Ranking\\
\hline
\\
\includegraphics[width=12mm, height=10mm]{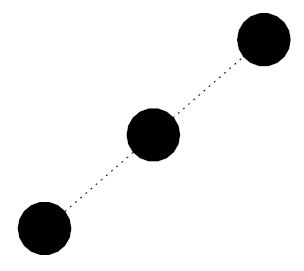} & (3, 2, 1) & (1, 2, 3)\\
\\
\includegraphics[width=12mm, height=10mm]{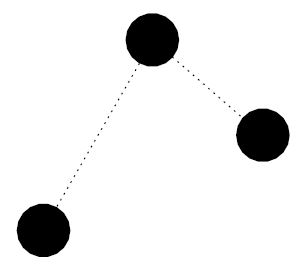} & (3, 1, 2) & (1, 3, 2)\\
\\
\includegraphics[width=12mm, height=10mm]{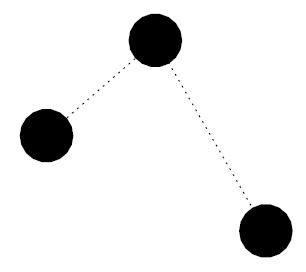} & (2, 1, 3) & (3, 1, 2)\\
\\
\includegraphics[width=12mm, height=10mm]{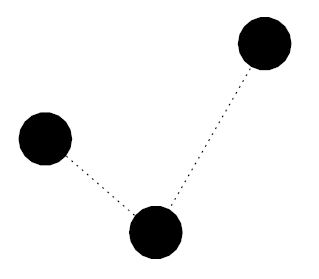} & (2, 3, 1) & (2, 1, 3)\\
\\
\includegraphics[width=12mm, height=10mm]{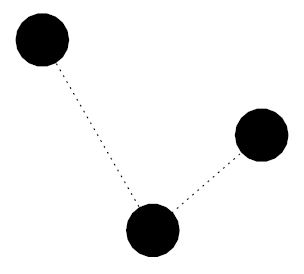} & (1, 3, 2) & (2, 3, 1)\\
\\
\includegraphics[width=12mm, height=10mm]{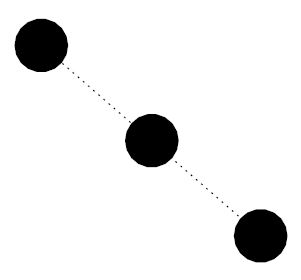} & (1, 2, 3) & (3, 2, 1)\\
\end{tabular}
\end{ruledtabular}
\end{table}

We can also consider non-overlapping points between these windows, denoted here by $w$.  That is, each $N$ point measurement time series is divided into  $\lfloor\frac{N-(m-1)\tau-1}{w} + 1\rfloor$ windows. Therefore, the situation with the fewest non-overlapping points ($w = 1$) results in the greatest number of windows, which is equal to $N - (m - 1)\tau$. When there is no overlap between the windows ($w = L = (m - 1)\tau$), the minimum number of windows is reached, equal to $\lfloor \frac{N}{(m - 1)\tau} \rfloor$.

\begin{figure*}[htbp]
\begin{tabular}{c}
\centerline{\includegraphics[width=0.85\textwidth]{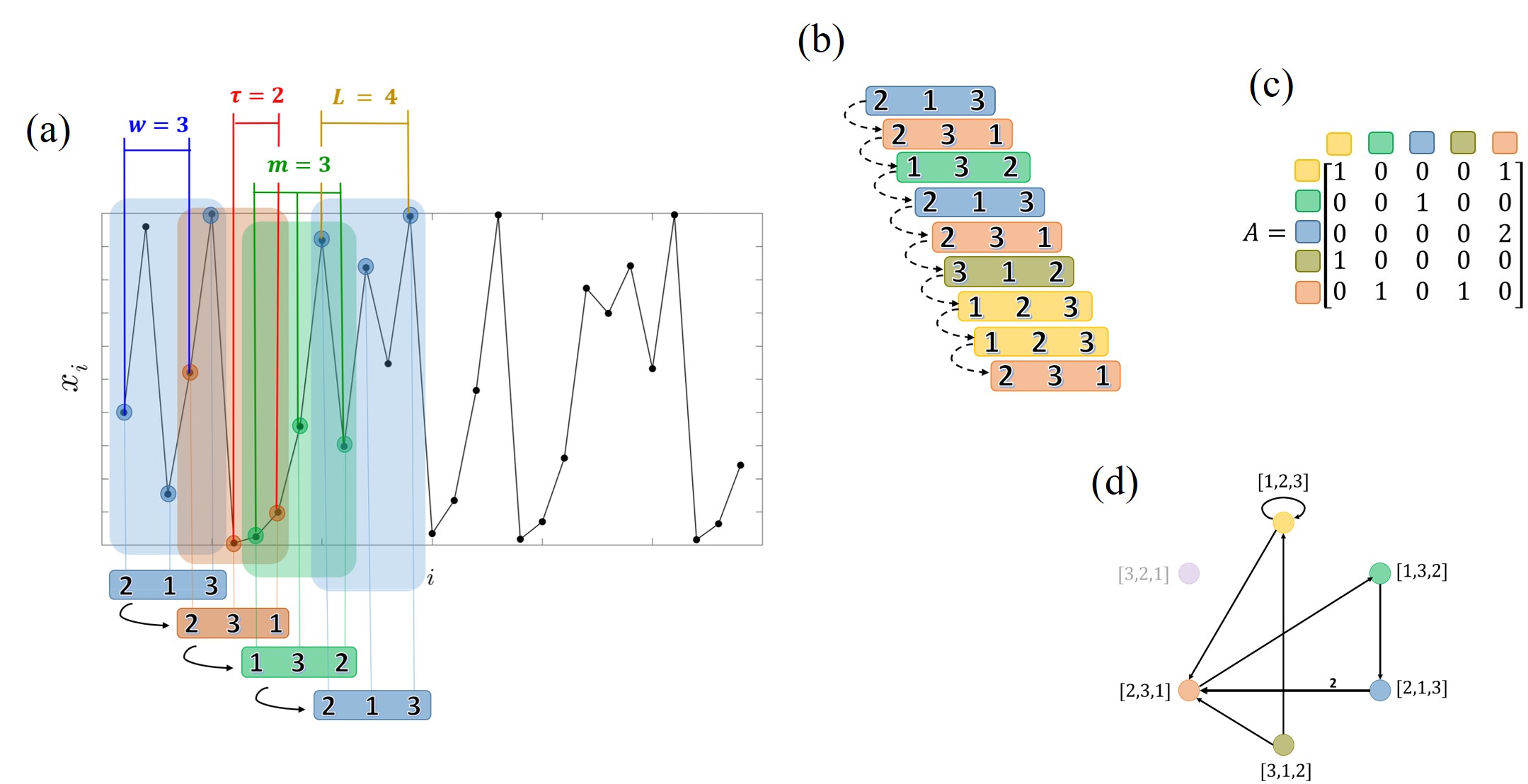}}
\caption{The process of constructing an OPN using a sample time series; a) finding the ordinal pattern by the embedding vectors with dimension $m$, lag $\tau$, the number of non-overlapping points of neighbor windows $w$, and the length of each window $L$; the ordinal patterns are found based on the chronological index rank of the selected samples; b) organizing ordinal symbols of the sequential windows; c) generating a network with adjacency matrix $A$ such that each of its nodes is a unique ordinal pattern, and the links are assigned between nodes based on the transitions between ordinal patterns in the time series; d) schematic of the constructed network. Since the example time series samples did not contain the ordinal pattern [3,2,1], the node is shown faded, and the final network does not contain a node for this ordinal pattern.}
\label{OPN}
\end{tabular}
\end{figure*}

Using chronological index ranking in each window, we can sort each window's selected $m$ points into at most $m!$ different orders as described above. Since some of the possible ordinal rankings may never occur within the time series, the network has at most $m!$ nodes. We can build a network for each time series in the following way. Associating each order of $z^{(i)}$ to a network node allows for considering a weighted directed graph. The nodes represent the states, and links between nodes illustrate the changes from each state to another. We also allow for self-loops when the state returns to itself. This network demonstrates how the ordered patterns in the sampled time series have changed over time. Figure~\ref{OPN} depicts this structure on a sampled time series of a map. The relative magnitude of $m$ points and their indices is taken into account in each window, and the adjacency matrix of the dynamical network is then created by arranging the sequential ordinal symbols.

To generate an OPN from a time series, we have the option of selecting three variables. The first is the window length, $L$, which determines the range of time series points to be analyzed. 
The second is the number of selected points in each window, $m$, which indicates the number of points in each window whose amplitudes we compare.
Obviously, considering more points in each window means having more nodes in the final network. 
It should only be considered that all of the parameters $L$, $m$, and $\tau$ are integers, and the values of $L$ and $m$ should be chosen in such a way that $\tau = \frac{L}{m-1}$ becomes an integer as well. 
The third is $w$, which represents the number of non-overlapping points of neighbor windows. In this study, we assumed that $w=1$, which means that the windows have the most overlap and that we require the most windows ($N - (m - 1)\tau$) to analyze the entire time series.

A reasonable prescription for selecting the most appropriate window length for extracting the dynamical features of the time series is to use traditional embedding methods. Because both tasks (reconstruction of phase space from the output time series, which is equivalent to a system's original phase space, and creating an OPN from a time series) practically revolve around the same central issue of determining how long the time series in a window can preserve the topological properties of the signal. In this regard, we use the method provided by Takens for embedding in this work \cite{r2}.

Although this is not the only possible solution, and in some cases, it does not work at all. For example, this method cannot be applied to delay systems with infinite dimensions. Also, in other non-delayed systems, which is challenging to find the appropriate dimension and delay for embedding the system in higher dimensions due to the complexity of their dynamics; any alternative method can be used for what window length is suitable for sequential partitions. Because of these characteristics, this method has fewer limitations than embedding methods and can be used to analyze systems with complex dynamics.

\subsection{\label{sec:First Return Maps}First Return Maps}

Since in the ordinal partitioning, each set of points is assigned to one and only one symbol, hence the ordinal partitions divide the system's attractor into several regions so that the union of all these regions covers the entire attractor. The border of each of these regions is defined by a set of points that separate these areas. We assert that each of these boundaries can be considered a Poincar\'e section on the attractor, and the FRM can be studied on it. Depending on which part of the attractor this ordinal partition cut is on, it can be considered as an appropriate or inappropriate Poincar\'e section on the attractor.

\subsection{\label{sec:Permutation Entropy}Permutation Entropy}

Permutation Entropy (PE) is a widely used measure of a dynamical system's complexity. It captures the permutation patterns and ordinal relationships between a time series' individual values. The ordinal patterns' probability distribution is then extracted. Since this criterion is non-parametric, it can be calculated directly from the one-dimensional time series and does not necessitate the use of a parametric model. It is also robust to noise, making it extremely effective at extracting the dynamic content of nonlinear time series \cite{r10}. By dividing the time series into a matrix of overlapping column vectors and generating the ordinal symbols (Fig.~\ref{OPN}a and b), the PE can be calculated. The $m$-dimensional vectors are thus mapped into distinct permutations of the ordinal rankings. Bandt and Pompe defined a time series' permutation entropy as the Shannon entropy of the corresponding set of ordinal symbols $s$ \cite{r11}: \ignore{\cite{r9}:}
\begin{equation}
h = - \sum_{i} p_i \log_2 p_i,
\label{eq.1}
\end{equation}
where $p_i = P (S = s_i)$ is the probability mass function for $s_i \in s$. Here, $s$ is the set of ordinal symbols in the symbolic dynamics $S$, and $s_i$ is equal to the relative occurrence of each symbol. To calculate the probability mass function, we consider that it is the stationary distribution of the Markov chain of $S$. So, each element of the matrix $p$ could be defined as:
\begin{equation}
p_{i,j} = \frac{a_{i,j}}{\sum_{k} a_{i,k}},
\label{eq.2}
\end{equation}
where $a_{i,j}$ are the number of transitions from state $i$ to state $j$ between at most $m!$ states that we can have in each window. Therefore, the probability mass function can be estimated in the following manner:
\begin{equation}
p_i \approx \frac{\sum_{j} a_{i,j}}{\sum_{k}\sum_{j} a_{k,j}}.
\label{eq.3} 
\end{equation}

To calculate the entropy for each set of ordinal partitions' points, we select all the points in the initial time series that have the same ordinal partition and combine them to form a new time series. The entropy for each ordinal partition is then calculated by comparing the amplitudes of the points in this new time series with each other and generating an OPN from it. We consider the window options $m' = 3$ and $\tau' = 1$ for creating the OPN from each ordinal partition's time series.
In this regard, for the first time, we proposed two new types of entropy that connect the probability of occurrence of each set of ordinal partition's points in the first time series and the OPN created from each of them. 

The first one is \textit{weighted entropy}, which uses the probability of occurrence of each set of ordinal symbols in the main time series as a scaling factor for $p_i$ as follows:
\begin{equation}
K = \frac{\mathcal{O}}{\mathcal{T}},
\label{eq.4_1} 
\end{equation}
where $\mathcal{O}$ shows the number of occurrences of each particular ordinal partition, and $\mathcal{T}$ stands for the total number of windows into which the time series is divided, which in this case, as was explained in Sec.\ref{sec:Ordinal Partition Networks}, is equal to $\mathcal{T} = N-(m-1)\tau$, where $N$ is the total length of time series and $m$ is the selected number of points for ordinal partitioning the main time series, and $\tau$ denotes the gap of each two consecutive selected points. By using this scaling factor, the weighted entropy is defined as follows:
\begin{equation}
\begin{split}
h_w & = -\sum_{i} K p_i \log_2 (K p_i) \\
 & = -K\sum_{i} p_i (\log_2 p_i + \log_2 K).
\end{split}
\label{eq.4} 
\end{equation}

The second entropy, which we call \textit{weighted transition entropy}, considers only the likelihood of entering a specific ordinal symbol among all available symbols as a scaling factor for $p_i$ as shown below:
\begin{equation}
\hat{K} = \frac{\hat{\mathcal{O}}}{\hat{\mathcal{T}}},
\label{eq.5_1} 
\end{equation}
where $\hat{\mathcal{O}}$ shows the number of \textit{entrances} to each particular ordinal partition, and $\hat{\mathcal{T}}$ denotes the total number of \textit{entering point} to all ordinal partitions. In this instance, for both $\hat{\mathcal{O}}$ and $\hat{\mathcal{T}}$, we only count the points that enter a new ordinal partition; if there are many sequential points in the same ordinal partition, only the first point is counted. The weighted transition entropy is defined using this scaling factor as follows:
\begin{equation}
\begin{split}
h_{wt} & = -\sum_{i} \hat{K} p_i \log_2 (\hat{K} p_i) \\
 & = -\hat{K}\sum_{i} p_i (\log_2 p_i + \log_2 \hat{K}).
\end{split}
\label{eq.5} 
\end{equation}

\section{\label{sec:Results and Discussions}Results and Discussions}

\begin{figure}[htbp]
\centerline{\includegraphics[width=\linewidth]{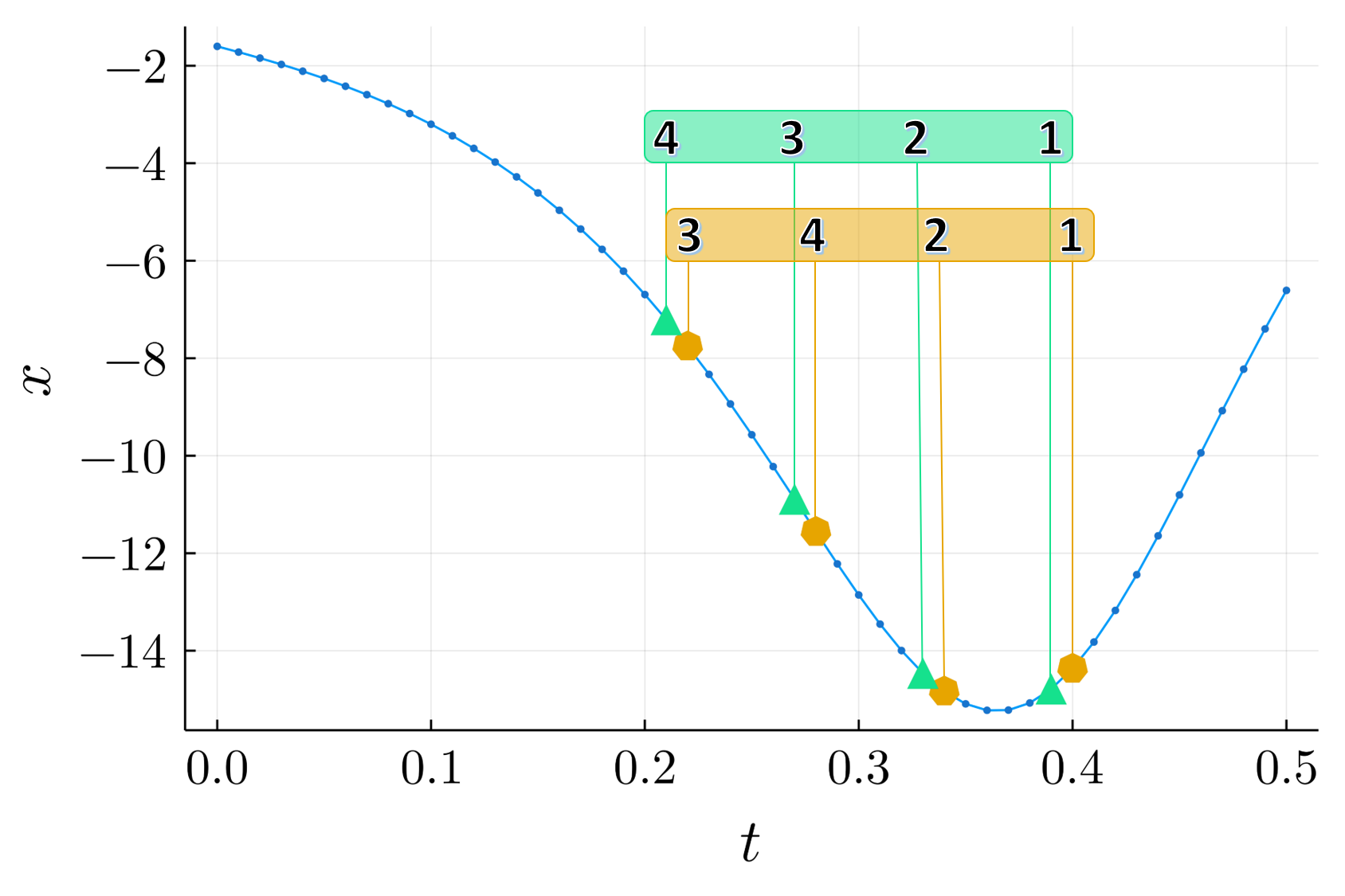}}
\caption{A small portion of the $x$ time series of the Lorenz attractor. The ordinal window is considered with a length of $L = 18$, $m = 4$ sampled points in each window, and $\tau=6$ points gap between each two. The non-overlapping points of neighbor windows are set to $w = 1$, which shows the ordinal partition window slides by one in each step. The points are sorted with chronological index ranking method.}
\label{LorenzTimeSeries_SampleOP}
\end{figure}

The $x$ time series of the Lorenz chaotic system has been chosen to check the efficacy of the presented method for time series analysis. The equations of the Lorenz system are as follows:
\begin{equation}
\label{Lorenz}
\begin{array}{lll}
\cfrac{dx}{dt} &=&\sigma (y - x), \\
\cfrac{dy}{dt} &=&x(\rho - z) - y, \\
\cfrac{dz}{dt} &=& xy - \beta z.
\end{array}
\end{equation}
where the parameters are set to $\sigma = 10$, $\rho = 28$, and $\beta = 8/3$ to ensure a chaotic regime. We calculated the system's response using random initial conditions for $10,000 \rm{sec}$ with a time resolution of $\Delta t = 0.01 \rm{sec}$. So, the time series response is obtained with $1,000,000$ points. We discard the first 90\% of this time series as a transition and focus on the last 10\%, which contains $100,000$ points. To study the system's dynamics, we first select the ordinal partition window length ($L$) based on Takens' embedding theorem window size, as explained in Sec.\ref{sec:Ordinal Partition Networks}. The number of sample points for each window ($m$) is determined by the computational constraints and the level of detail that we want to extract from the system. To accomplish this, we choose the embedding dimensions $M = 3$ and the embedding lag $T = 9$ points, which correspond to roughly one-fourth of the Lorenz $x$ time series period according to the considered time resolution. As a result, we suppose that each ordinal window has a length of $L = (M-1)T = 18$, and we assume the number of sample points for each window $m = 4$, which gives us the ordinal lag $\tau = \frac{L}{m-1} = 6$ for each window. Figure~\ref{LorenzTimeSeries_SampleOP} shows a small part of the Lorenz time series and two successive ordinal partition windows on it. After mapping each window to a symbol, the result will be assigned to the first node of the window on the time series.

\subsection{\label{sec:Weighted Entropy and Weighted Transition Entropy}Weighted Entropy and Weighted Transition Entropy}

\begin{figure}[htbp]
\centerline{\includegraphics[width=\linewidth]{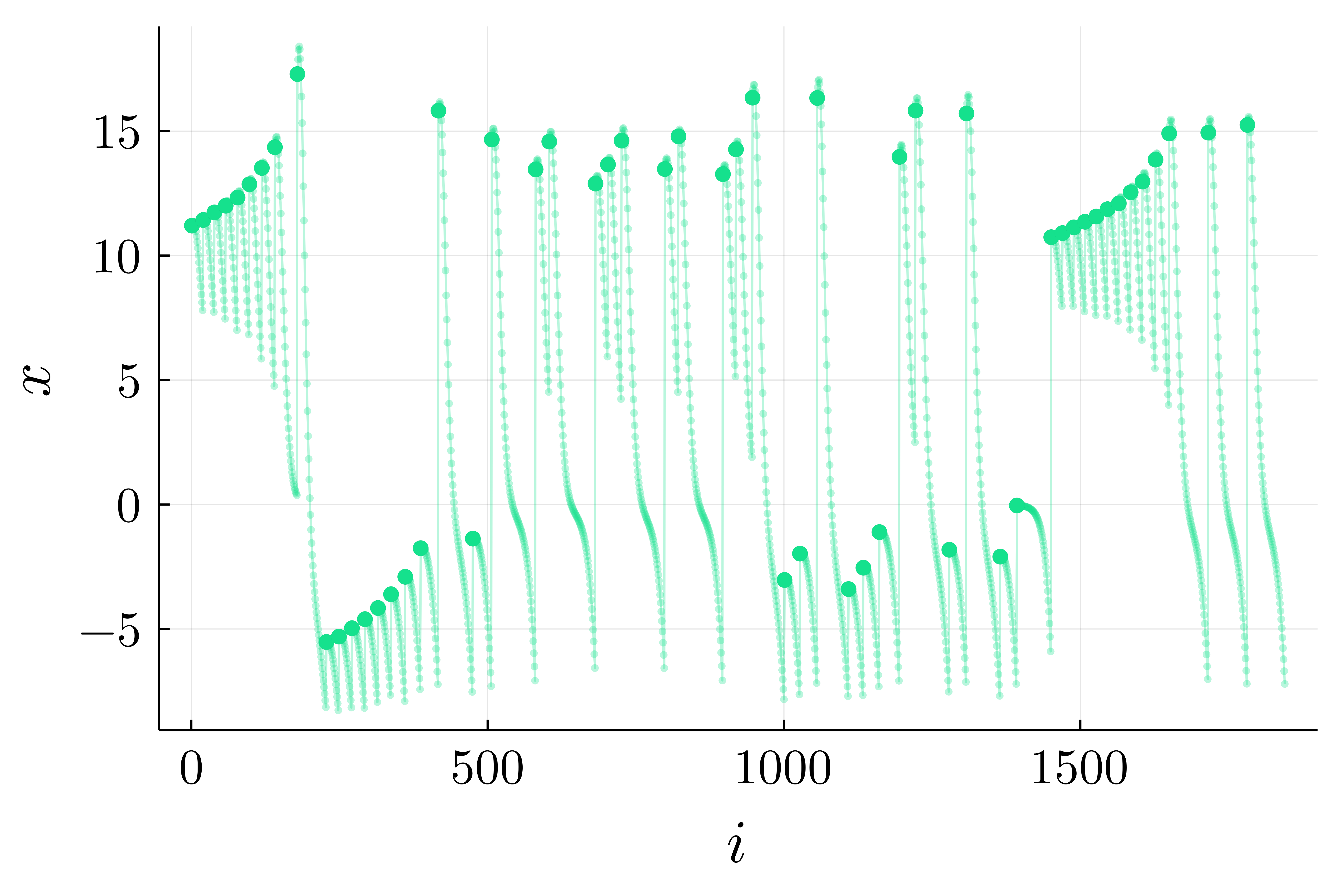}}
\caption{Occurrence of the ordinal partition [4,3,2,1] derived from the Lorenz system's $x$ time series. The larger dots represent the entry points to this ordinal partition from other ordinal partitions.}
\label{a_LorenzTimeSeries_descending}
\end{figure}

\begin{figure*}[htbp]
\begin{tabular}{cc}
\begin{tabular}{c}
        {\includegraphics[width = 0.45\linewidth]{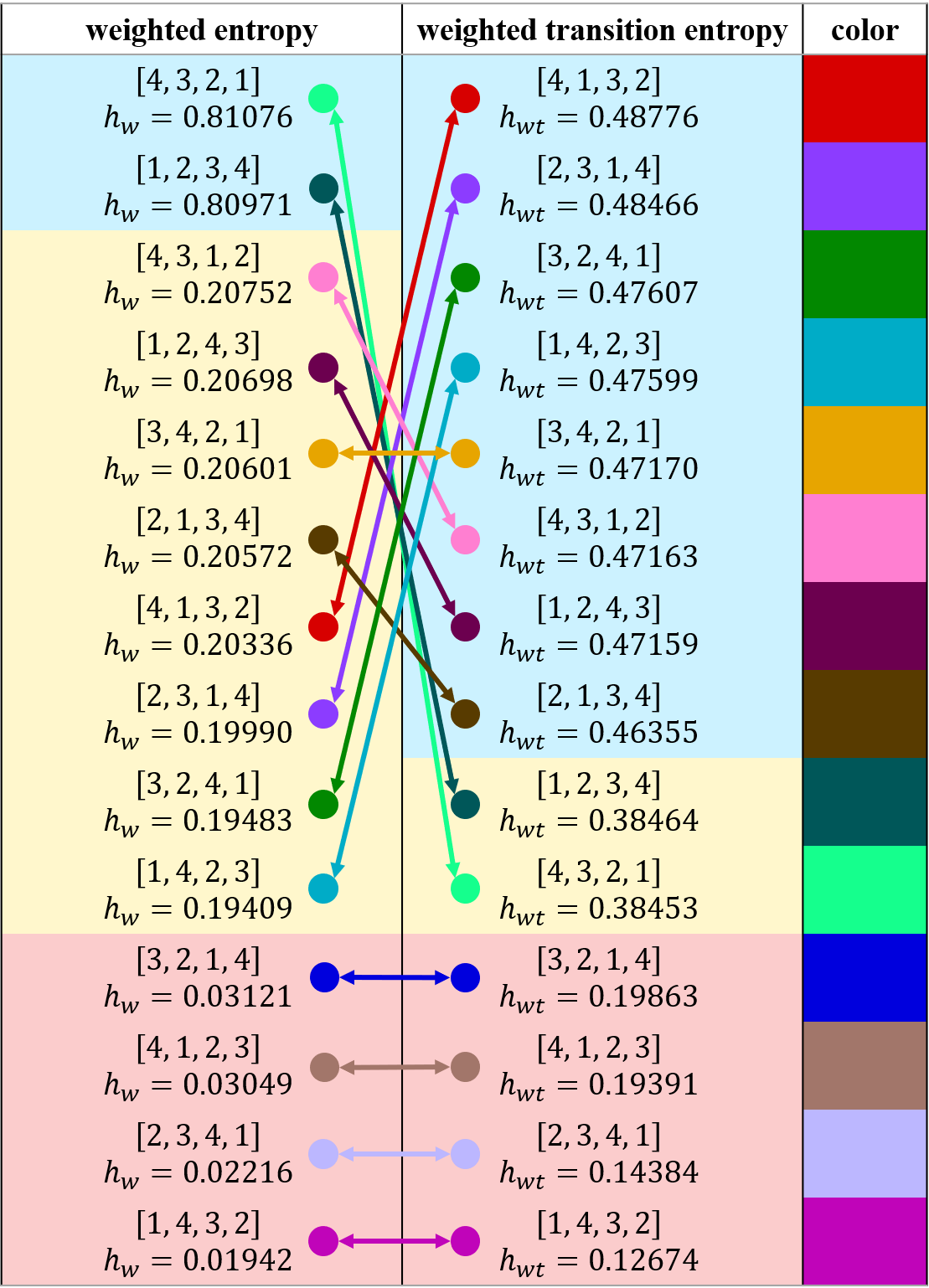}} \\ (a) 
\end{tabular} 
&
\begin{tabular}{c} 
        {\includegraphics[width = 0.45\linewidth]{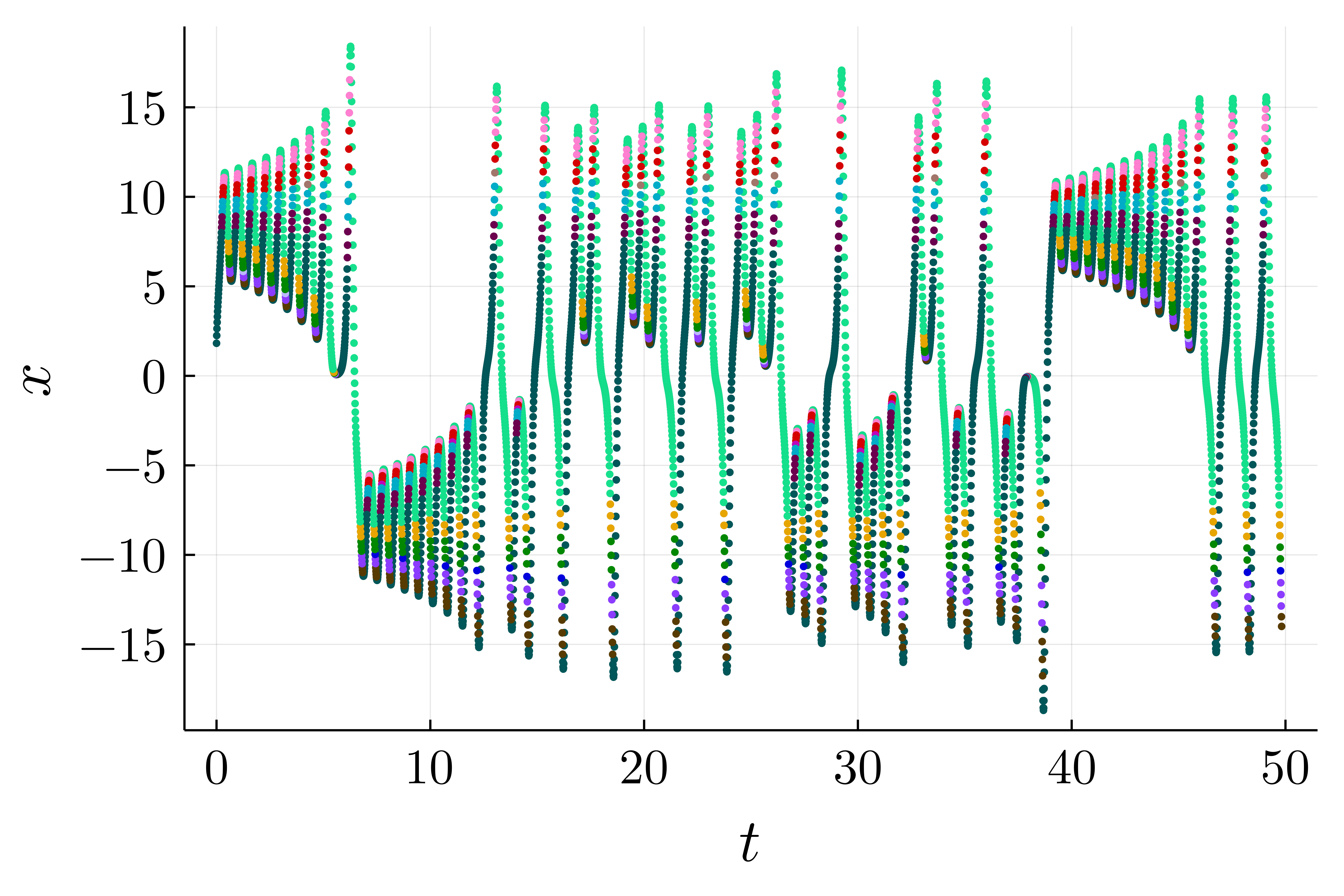}}\\ (b)\\
        {\includegraphics[width = 0.45\linewidth]{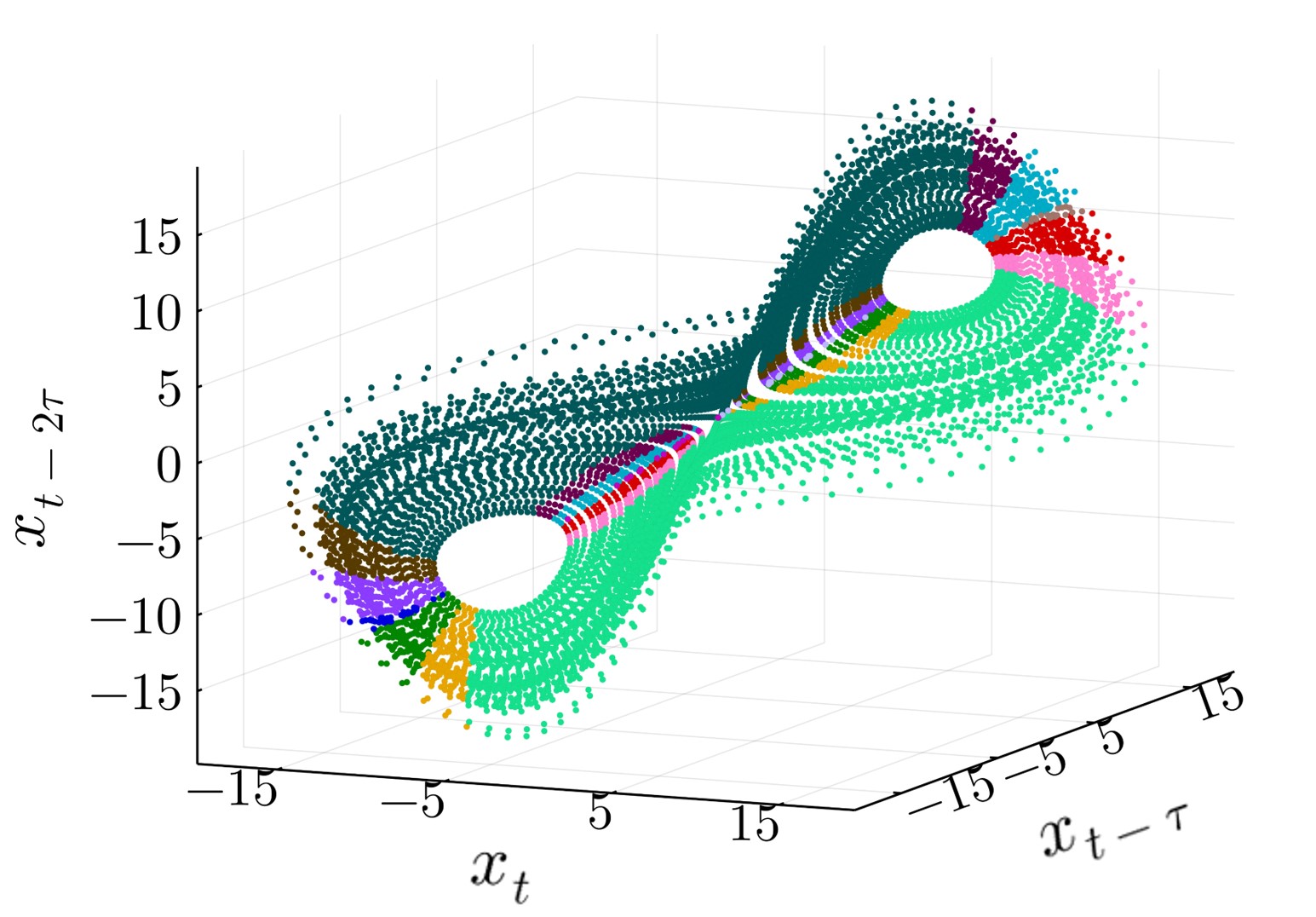}}\\ (c)
\end{tabular}
\end{tabular}

    \caption{(a) The list of ordinal partitions that have occurred on the Lorenz attractor. The first column is ordered by weighted entropy ($h_w$), and the second column is ordered by weighted transition entropy ($h_{wt}$). Each ordinal partition has been assigned a specific color, as shown in the third column. The background colors in the first and second columns represent the entropy level in their respective groups (See Fig.~\ref{SortedEntropy}). Blue cells are in the first level of entropy, yellow cells are in the second level, and red cells are in the third level. (b) Lorenz's $x$ time series colored based on the ordinal partition. Each window's ordinal partition has been assigned to the first point of the window. (c) The Lorenz embedded attractor from its $x$-colored time series using Takens' embedding theorem with an embedding dimension of $M=3$ and an embedding lag of $T = 9$. The high entropy sections serve as appropriate Poincar\'e sections on the attractor.}
    \label{LorenzColoredOP}
\end{figure*}

To divide the time series into different ordinal partitions, we consider a window as shown in Fig.~\ref{LorenzTimeSeries_SampleOP} and move it through the entire time series. This allows us to assign a symbol to each time series point. We then generate an OPN from each of these symbols to determine the $h_w$ and $h_{wt}$ of each collection. 
Figure~\ref{a_LorenzTimeSeries_descending} depicts the descending node's time series, which has been separated from Lorenz's time series as an independent new time series. Comparing Fig.~\ref{a_LorenzTimeSeries_descending} and Fig.~\ref{LorenzColoredOP}b shows that this ordinal of the time series is correctly limited in the range of [-8.3 18.5]. To compute $h_w$ and $h_{wt}$ for this new time series, we first create an OPN with the parameters $m' = 3$ and $\tau' = 1$ out of it. Then, to calculate the scaling factors for this specific ordinal partition, we set $\mathcal{O}$ equal to the length of this time series segment and set $\hat{\mathcal{O}}$ equal to the number of large points in it that specify the first entries to this section. As previously stated, $\mathcal{T}$ and $\hat{\mathcal{T}}$ are calculated from the entire time series of the system and are the same for all ordinal partitions, not just the partition under study.

Figure~\ref{LorenzColoredOP}a displays the $h_w$ for each ordinal partition, which is sorted by highest to lowest entropy in the first column. The $h_{wt}$ is displayed similarly in the second column. A color has been assigned to each ordinal partition, which is more clearly shown in the third column. This allows us to assign a color to each time series point, as seen in Fig.~\ref{LorenzColoredOP}b. With this colored time series, we may recreate the attractor in three dimensions according to Takens' embedding theorem, yielding Fig.~\ref{LorenzColoredOP}c. By comparing the topological properties of this attractor to those of the Lorenz original attractor in Fig. 1a,
it can be seen that the embedded attractor mimics the primary features of Lorenz, such as its two wings and two holes. As can be seen, each ordinal partition has a distinct location on the attractor; thus, the entry points of ordinal partitions with high entropy values can be thought of as Poincar\'e sections on the attractor. On the other hand, the ordinal partitions with low entropy values are just a few single dots on the time series, randomly arranged on the embedded attractor, and cannot be considered as a good section that captures the dynamic of the system.

Attractor regioning by ordinal partitions occurs in all dynamical systems, regardless of their dimension or complexity. The appendix shows similar results for the R{\"o}ssler and Mackey-Glass systems as examples. It is obvious that because Mackey-Glass has a more complicated dynamic, it experiences more different ordinal partitions.

\begin{figure}[htbp]
\centerline{\includegraphics[width=\linewidth]{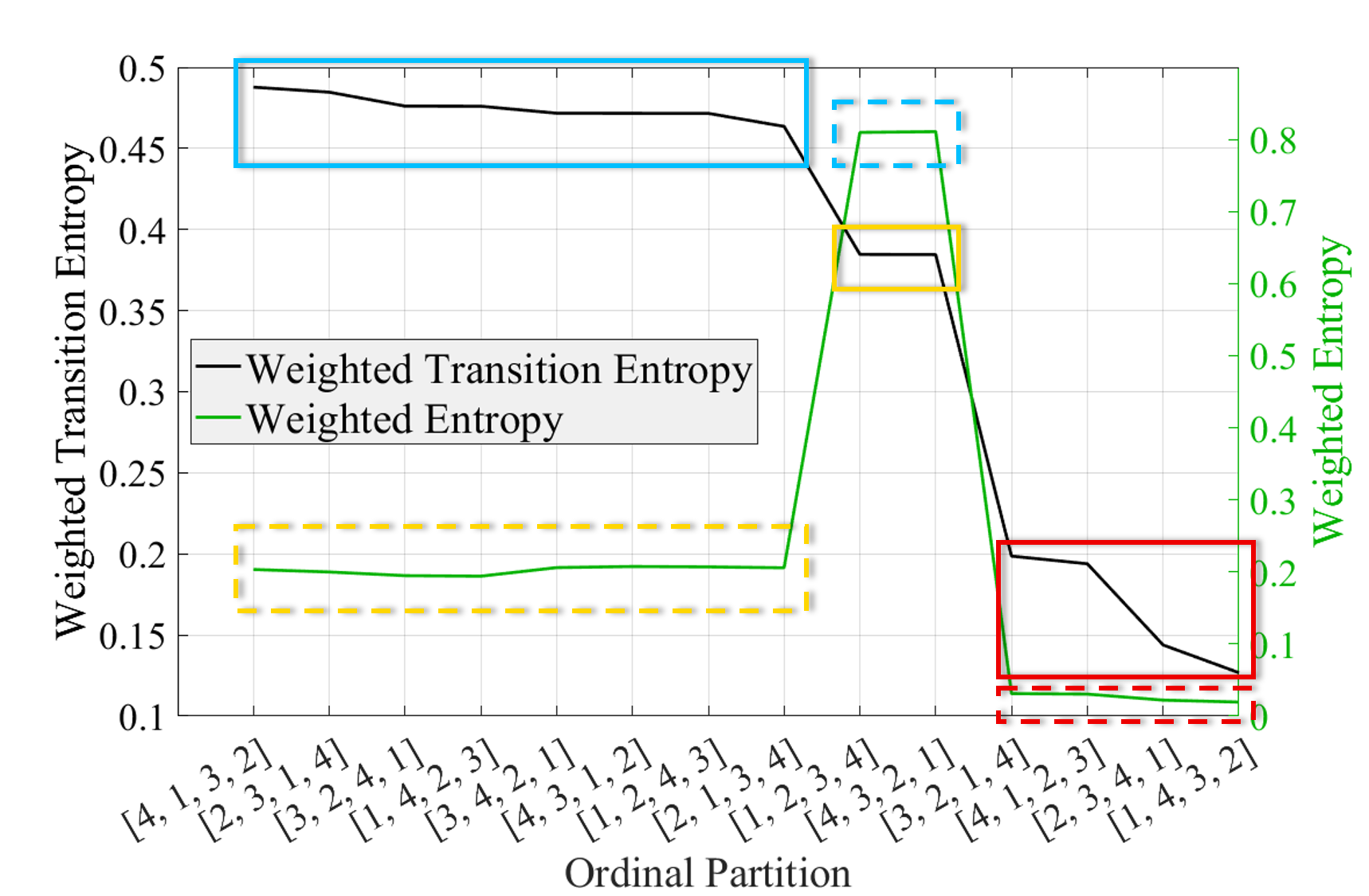}}
\caption{The weighted transition entropy of each ordinal partition sorted in descending order (black line) and their corresponding weighted entropy (green line). The solid boxes represent different levels of weighted transition entropy, and the dashed boxes represent different levels of weighted entropy: blue for the first level, yellow for the second level, and red for the third level of entropy. The diagram shows that the first and second levels of entropy are switched in weighted transition entropy and weighted entropy, but the third level remains constant.}
\label{SortedEntropy}
\end{figure}

Figure~\ref{SortedEntropy} depicts the amount of $h_{wt}$ and $h_w$ of each ordinal partition, which are sorted in descending order based on $h_{wt}$. The range of $h_{wt}$ is shown on the left axis, and the range of $h_w$ is shown on the right axis. Although the ranges of these two entropies differ, both show three levels of entropy. The $h_{wt}$ is represented by the black line, and the $h_w$ is represented by the green line. The levels of $h_{wt}$ are shown in solid boxes, and the levels of $h_w$ are shown in dashed boxes. In both, blue boxes represent the first level of entropy, yellow boxes represent the second level, and red boxes represent the third level. The graph illustrates that the first and second entropy levels are switched in different entropy methods, but the third level has the same group of ordinal partitions. At the first level, $h_w$ has a higher value than $h_{wt}$, but $h_{wt}$ has a higher value at the second and third levels. Generally, the variance of $h_w$ is greater than the variance of $h_{wt}$. This is because the weight of the ascending and descending ordinal partitions (jade green and light green parts in Fig.~\ref{LorenzColoredOP}c) in $h_w$, which is calculated by the number of all points that have that ordinal, is much higher than the weight of all other ordinal partitions. However, the $h_{wt}$ calculates the weights only based on the entry points to each section, putting them in the same range.

As shown in the appendix, the number of levels and their entropy amplitude in both $h_w$ and $h_{wt}$ can be changed due to the system's dynamic and the number of points considered in each window of the ordinal partition. The changes in $h_{wt}$ can be smoother as more points are considered in the ordinal window, allowing the first and second levels of entropy to be merged with each other \ignore{(Fig.~\ref{LorenzColoredOP_m10:a})}. Furthermore, depending on the system's dynamic, the number of levels in $h_w$ and $h_{wt}$ may differ \ignore{(Fig.~\ref{RosslerColoredOP:a} and Fig.~\ref{RosslerColoredOP:b})}, or there may be no obvious levels at all\ignore{(Fig.~\ref{MackeyGlassColoredOP:a} and Fig.~\ref{MackeyGlassColoredOP:b})}.

\subsection{\label{sec:First Return Maps According to Entropy}First Return Maps According to Entropy}

\begin{figure}[htbp]
    \centering
    \sidesubfloat[]
        {\includegraphics[width = 0.9\linewidth]{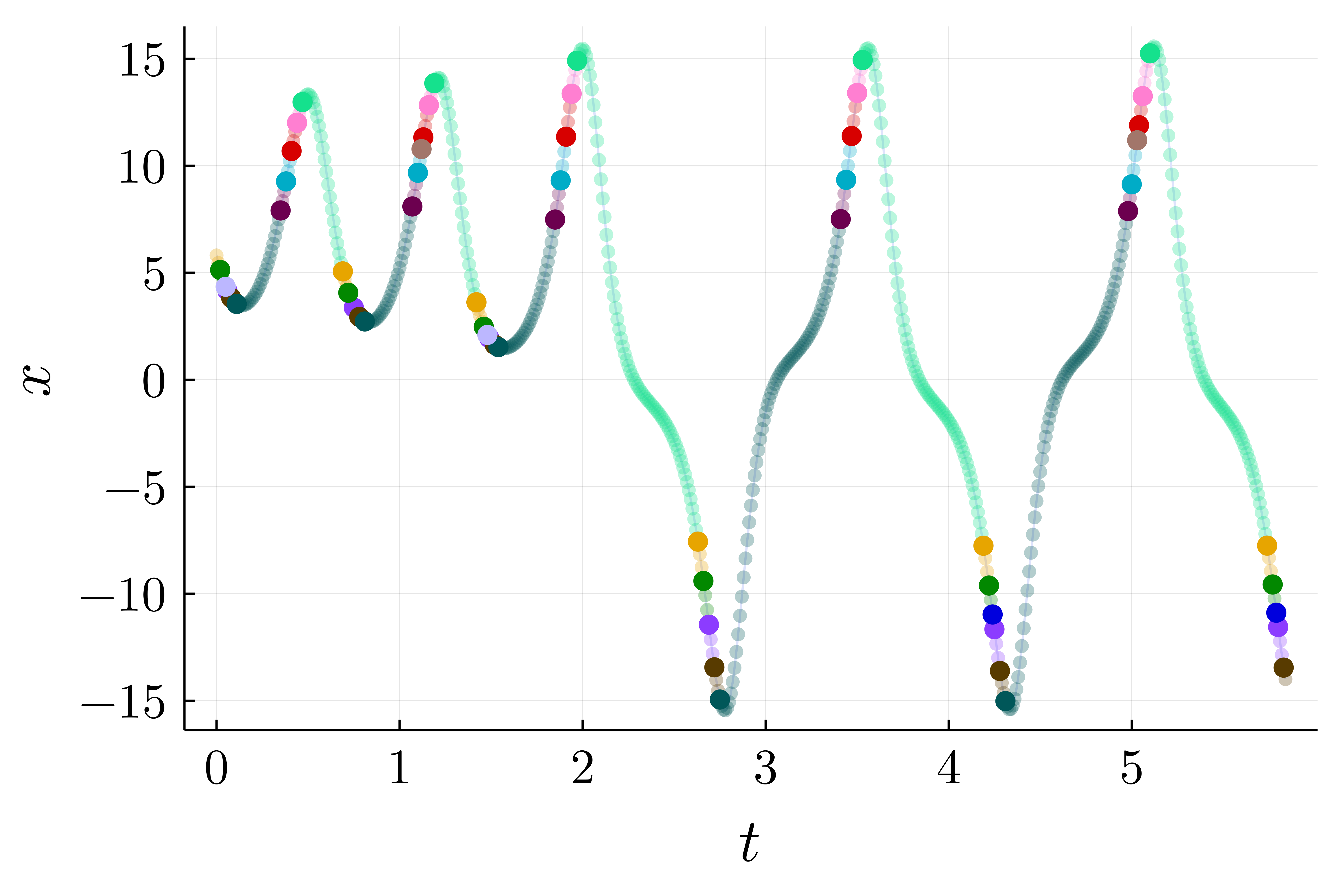}
        \label{LorenzColoredTransparentOP:a}}\\[-1pt]
   \sidesubfloat[]
        {\includegraphics[width = 0.9\linewidth]{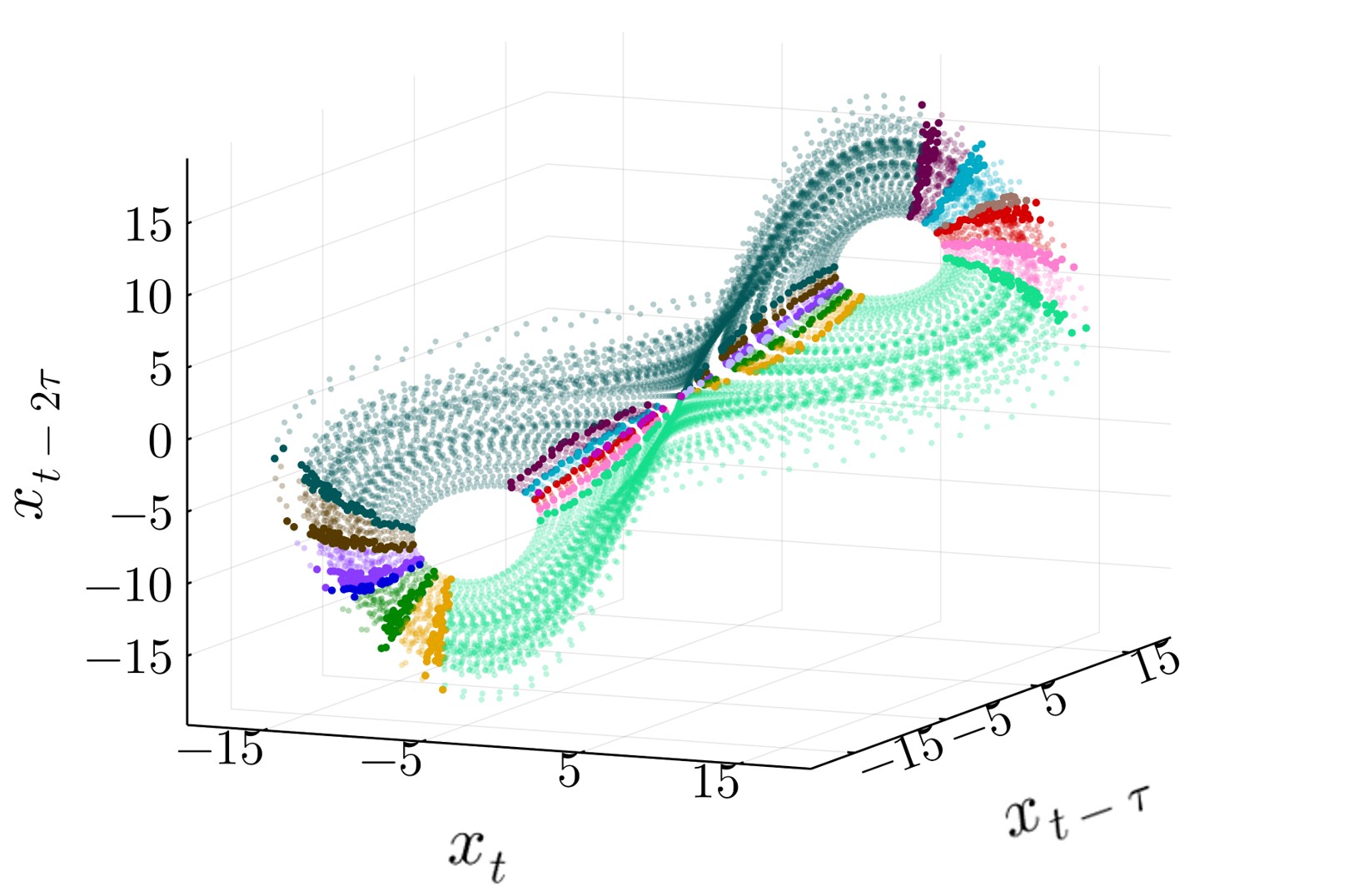}
        \label{LorenzColoredTransparentOP:b}}\\[-1pt]
    \sidesubfloat[]
        {\includegraphics[width = 0.9\linewidth]{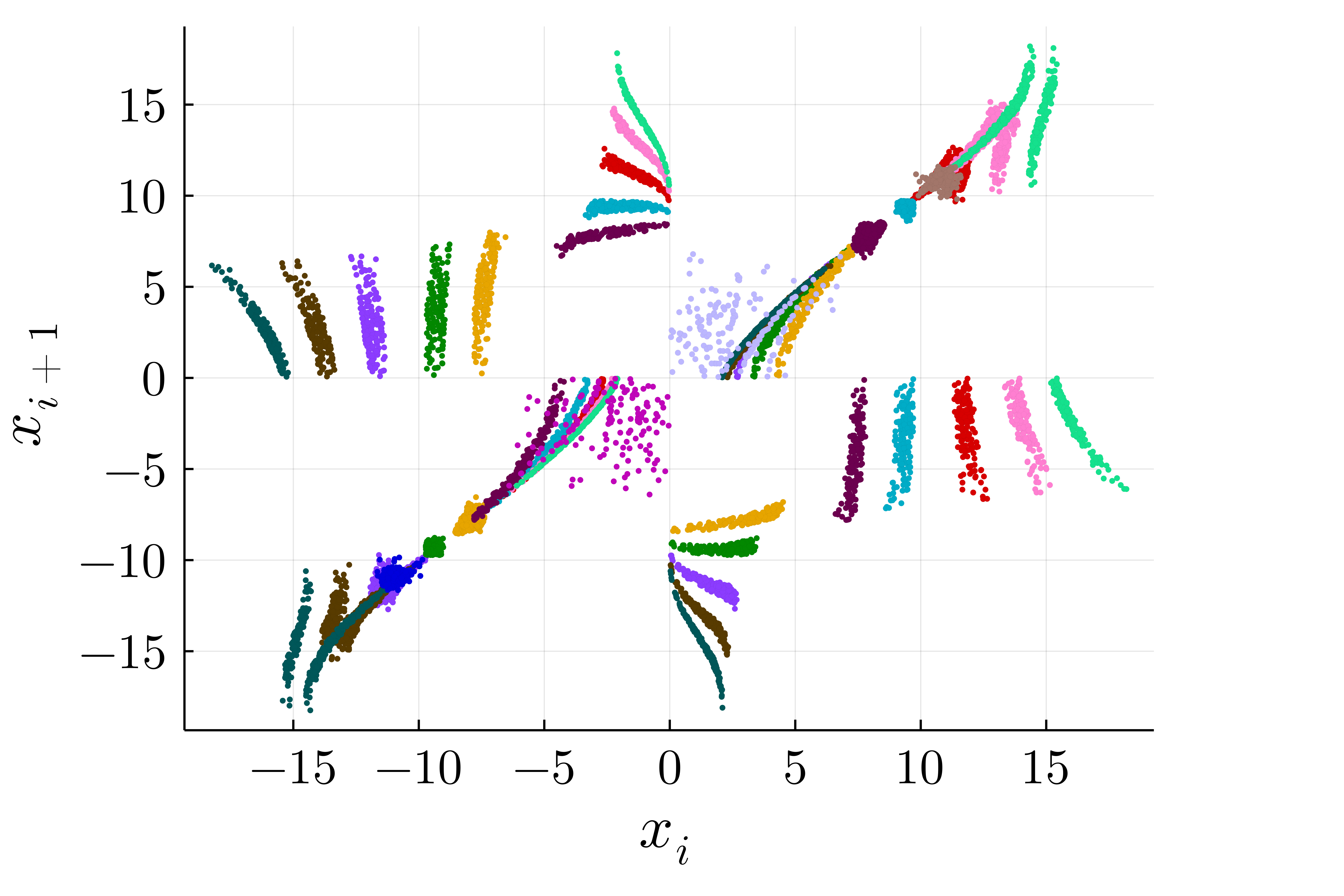}
        \label{LorenzColoredTransparentOP:c}}
    \caption{(a) Lorenz's $x$ weighted transition colored time series based on the ordinal partition. Larger dots represent the first points on the time series that enter an ordinal partition. The remaining points in the same ordinal partition are plotted with smaller, transparent dots. (b) The Lorenz embedded attractor from Fig.~\ref{LorenzColoredTransparentOP:a} time series using Takens' embedding theorem with an embedding dimension of $M=3$ and an embedding lag of $T = 9$. The entry points to each section (which are used to calculate weight in weighted transition entropy) are plotted larger and non-transparently. (c) Corresponding FRMs of each ordinal partition. N.B. Figure.~\ref{LorenzColoredTransparentOP:b} is not required to draw Fig.~\ref{LorenzColoredTransparentOP:c}, and it is only provided to demonstrate why high-entropy ordinal partitions can be considered good Poincar\'e sections.}
    \label{LorenzColoredTransparentOP}
\end{figure}

Figure~\ref{LorenzColoredTransparentOP} shows the entry points for each section more clearly, with larger and non-transparent points both on the time series (Fig.~\ref{LorenzColoredTransparentOP:a}) and the embedded attractor (Fig.~\ref{LorenzColoredTransparentOP:b}). These points are the ones that are used to calculate the weight in $h_{wt}$. Each set of dots in each ordinal partition has created an FRM, which is shown in Fig.~\ref{LorenzColoredTransparentOP:c}. Comparing Fig.~\ref{LorenzColoredTransparentOP:c} and Fig.~\ref{LorenzPoincareSection:c} indicates that the FRM of the descending ordinal partition [4, 3, 2, 1] (light green) gives the same FRM as the local maxima. Furthermore, the other FRMs obtained from other ordinal partitions in the first and second entropy levels are also precise for representing the main features of the attractor. While the FRMs of the third entropy level are just some random dots, and they don't represent any particular information about the topological structure of the attractor.

As previously stated, the attractor of any dynamical system can be divided into different ordinal partitions based on any desired number of points in the ordinal window. As a result, the FRMs derived from these ordinal partitions can represent the attractor's dynamical features. The difference between the FRMs of different ordinal partitions grows larger as the attractor's dynamic becomes more complicated. 
As shown in the appendix, the difference between the FRMs of different ordinal partitions for R{\"o}ssler is greater than the difference between the FRMs of different ordinal partitions for Lorenz and this difference is greater in the Mackay-Glass system compared to both (as it has a much more complicated dynamics).

\begin{figure}[htbp]
    \centering
    \sidesubfloat[]
        {\includegraphics[width = 0.9\linewidth]{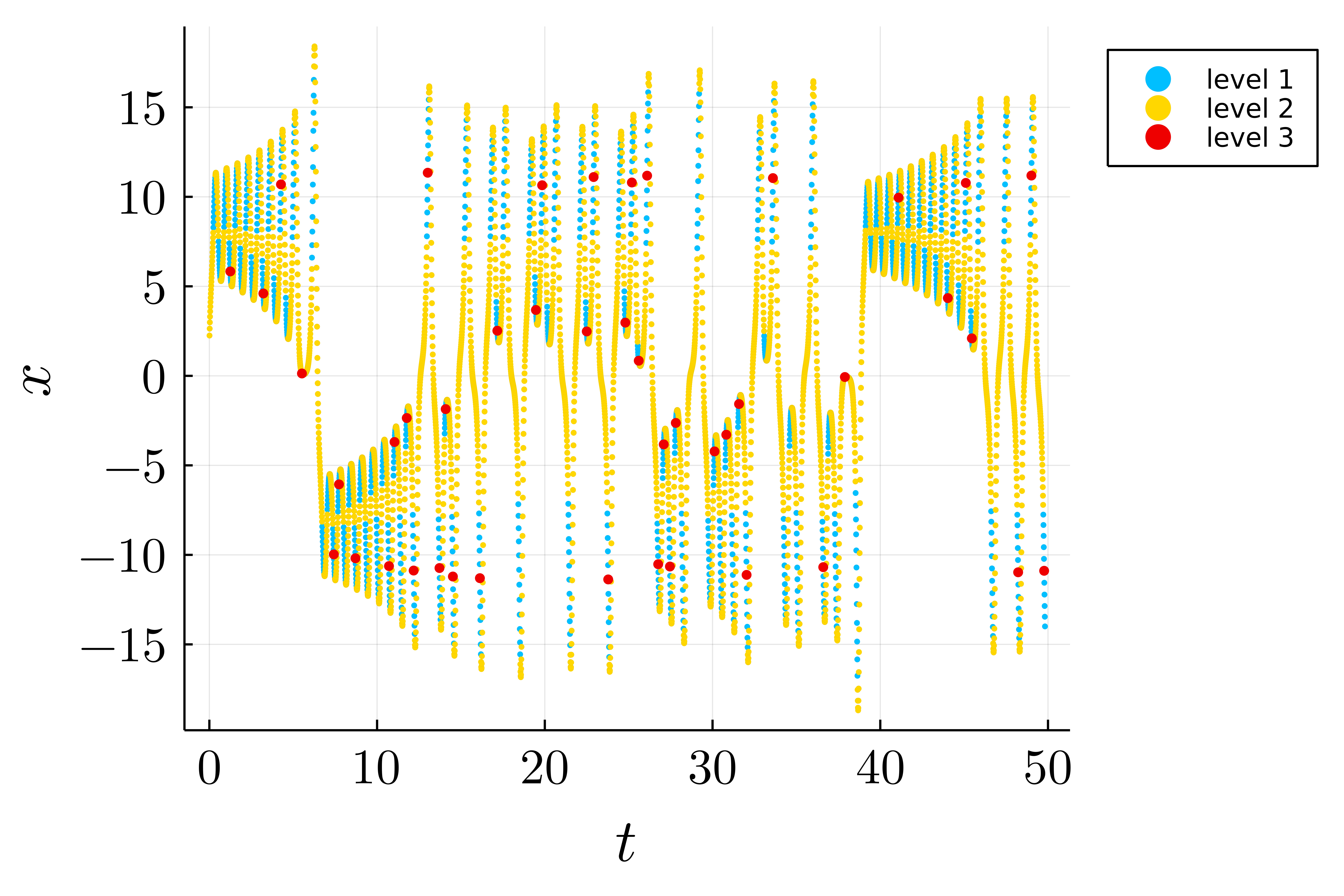}
        \label{LorenzColoredLevel:a}}\\[-1pt]
    \sidesubfloat[]
        {\includegraphics[width = 0.9\linewidth]{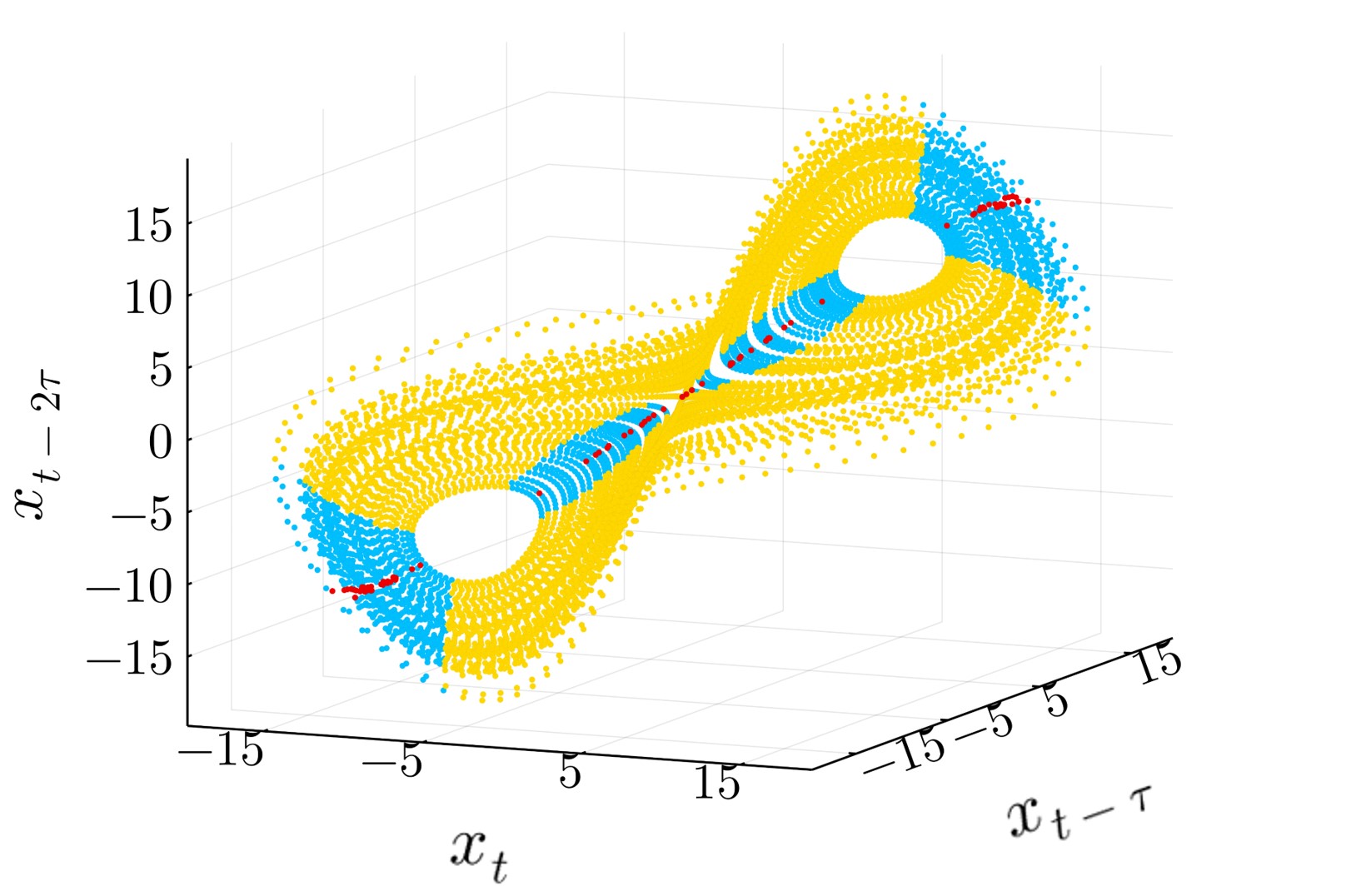}
        \label{LorenzColoredLevel:b}}\\[-1pt]
   \sidesubfloat[]
        {\includegraphics[width = 0.9\linewidth]{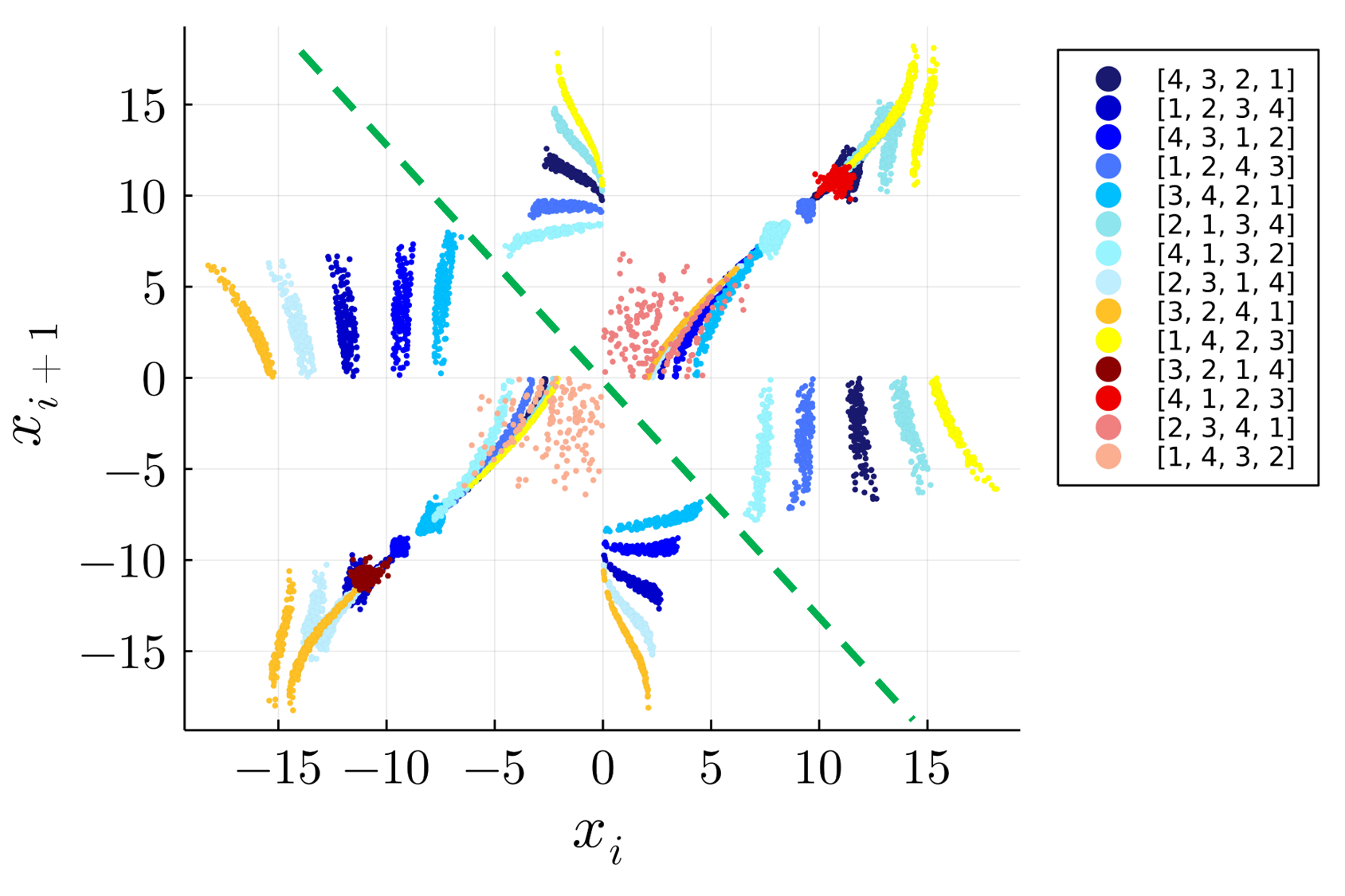}
        \label{LorenzColoredLevel:c}}\\[-1pt]
    \caption{(a) Lorenz's $x$ time series colored according to the level of weighted transition entropy of each ordinal partition. Blue-shaded dots represent the first level of entropy, yellow-shaded dots represent the second level, and red-shaded dots represent the third level. The red dots clearly include a small number of points on the time series. (b)  The Lorenz embedded attractor from Fig.~\ref{LorenzColoredLevel:a} time series using Takens' embedding theorem with an embedding dimension of $M=3$ and an embedding lag of $T = 9$. (c) Corresponding FRMs of each ordinal partition colored by the level of weighted transition entropy. The first-level FRMs are plotted in blues, from dark blue (the highest entropy value in the first level) to light blue (the lowest entropy in the first level). Similarly, the FRMs of the second and third entropy levels have been plotted in yellows and reds, respectively. The green diagonal dashed line divides the FRM's surface into two parts, each of which corresponds to one wing of the attractor.}
    \label{LorenzColoredLevel}
\end{figure}

Figure~\ref{LorenzColoredLevel} more clearly depicts the position of each level of entropy on the time series (Fig.~\ref{LorenzColoredLevel:a}) and on the embedded attractor (Fig.~\ref{LorenzColoredLevel:b}). It can be seen that the level 1 (blue) and level 2 (yellow) points cover almost the entire time series and the attractor, and level 3 (red) points are only a few dots between them. The FRMs of all ordinals in each of these levels are shown in Fig.~\ref{LorenzColoredLevel:c}. The FRMs in each level are colored with a color spectrum ranging from dark to light, from the highest entropy in the level to the lowest entropy. The FRMs are divided into two wings of the attractor by the green diagonal dashed line in the middle of Fig.~\ref{LorenzColoredLevel:c}. It can be seen that all of the FRMs on the first and second levels have a section above and a section below this diagonal. However, each group of red dots belonging to an ordinal partition of the third level of entropy only has points on one side of the diagonal, indicating that their relative sections on the attractor cannot be used as suitable Poincar\'e sections to study the system.

\begin{figure}[htbp]
\centerline{\includegraphics[width=\linewidth]{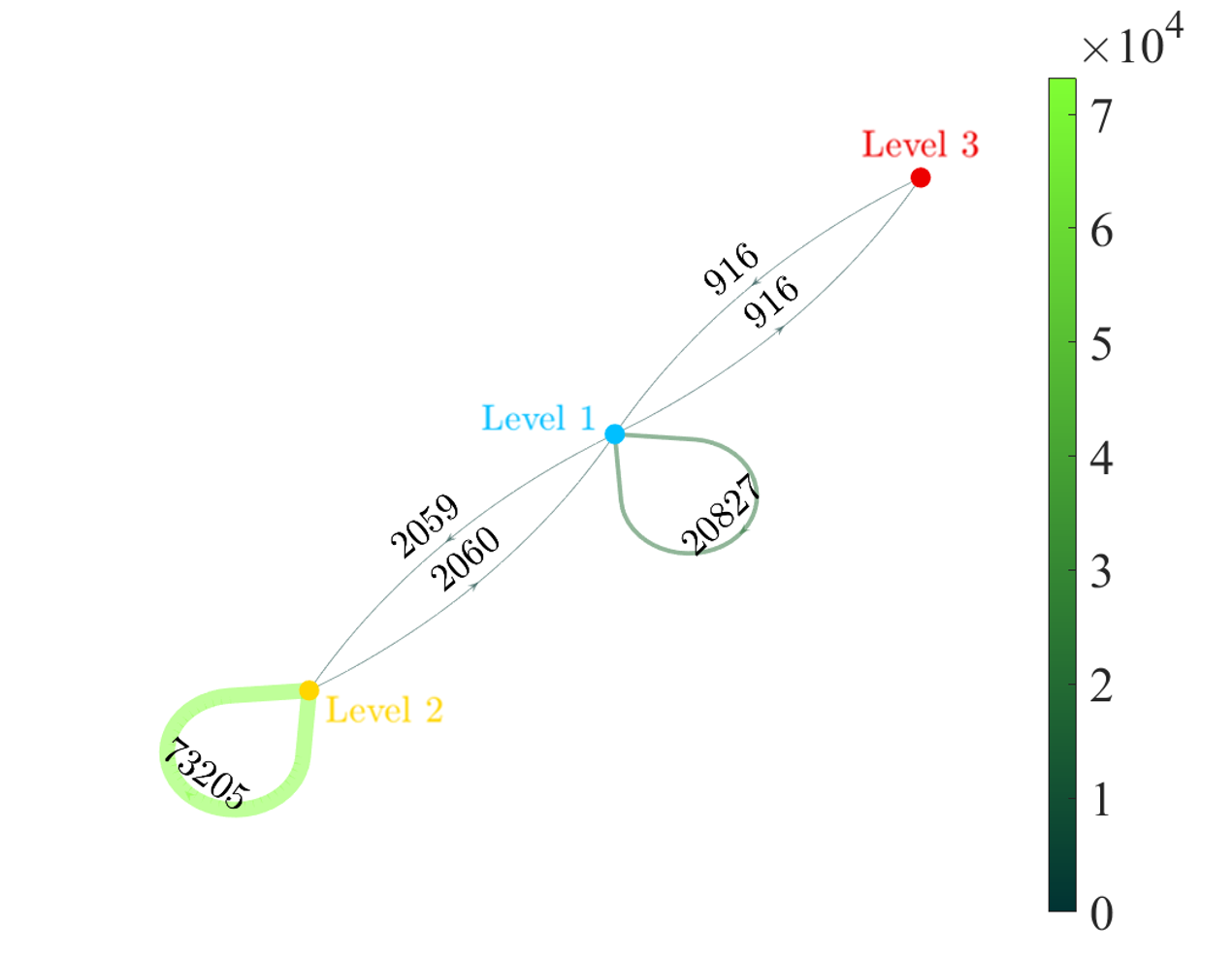}}
\caption{Network of weighted transition entropy levels. The nodes represent the entropy levels, and the links represent the number of transitions between every two levels or from one level to itself. Each link's thickness and color correspond to its weight.}
\label{NetworkLevels}
\end{figure}

With greater precision in Figs.~\ref{LorenzColoredLevel:a} and \ref{LorenzColoredLevel:b}, it can be understood that the points of level 3 are always placed between the points of level 1. This means that it never moves directly from level 2 to level 3. This fact inspires the idea of constructing a weighted directed network from the system, where each node represents a level of entropy, and the links represent the transition between entropy levels over time. Figure~\ref{NetworkLevels} shows this network which is constructed using $h_{wt}$ levels. It clearly shows that there is no direct path from level 2 to level 3, and there are no self-loops on the third level node. Level 2 has the most self-loops on itself, which makes sense when looking at Fig.~\ref{LorenzColoredLevel:b}. The attractor is mostly covered by yellow dots, representing the second entropy level.

\section{\label{sec:Conclusion}Conclusion}

We presented a new generic method for analyzing dynamical systems based on their output time series that can be applied to a wide range of dynamical systems regardless of complexity. This method extracts Poincar\'e sections by using ordinal partitions without the use of an embedded attractor or generating a high-dimensional complex network.
The novel types of weighted permutation entropies presented here help to rank the effectiveness of these Poincar\'e sections. Since ordinal partitions with high entropy values have a larger sector on the attractor, they can be considered suitable Poincar\'e sections for analyzing the system's dynamical behaviors.
The combination of the FRMs from all of these sections can provide a group of system information that aids in analysis. 
A group of FRMs can either allow one to choose their favorite or provide a more comprehensive picture of the dynamics.

We investigate the application of this method to various dynamical systems and, in all cases, are able to represent several one-dimensional maps of the system accurately, each capturing the system's main dynamical characteristic.
In particular, our method generates various types of FRMs that show detailed aspects of the system's dynamic in more complicated dynamical systems (such as Mackey-Glass). Embedding demonstrates why all of the high entropy ones are good and capture the system accurately. Some of these good FRMs are comparable to the well-known Poincar\'e sections (like $\dot{x}=0$).

While the method we have proposed here is directed at the classic construction of FRMs --- successive intersections of Poincar\'e sections --- the method itself could be applied without modification to any time series, making our approach superior to traditional methods such as finding maxima. 
Without the underlying continuity of a flow, we would only be able to look at successive first inclusion within a particular ordinal class. This would still represent a projection to a lower dimensional system, but the classical connection between flow and the first return map has no obvious immediate analog. This may also yield a valuable technique for discontinuous time series analysis, but the analytic underpinnings are less obvious (to us).


\begin{acknowledgments}
We wish to acknowledge Australian Research Council funding for DP200102961. SDA also acknowledges funding from the Forrest Foundation.
\end{acknowledgments}

\appendix

\section{Results on other dynamical systems}

This section provides examples of our method's application to other dynamical systems.

\subsection{\label{app:Lorenz_m10}Lorenz with $m=10$}

\begin{figure*}[htbp]
\begin{tabular}{cc}
\begin{tabular}{c}
        {\includegraphics[width = 0.49\linewidth]{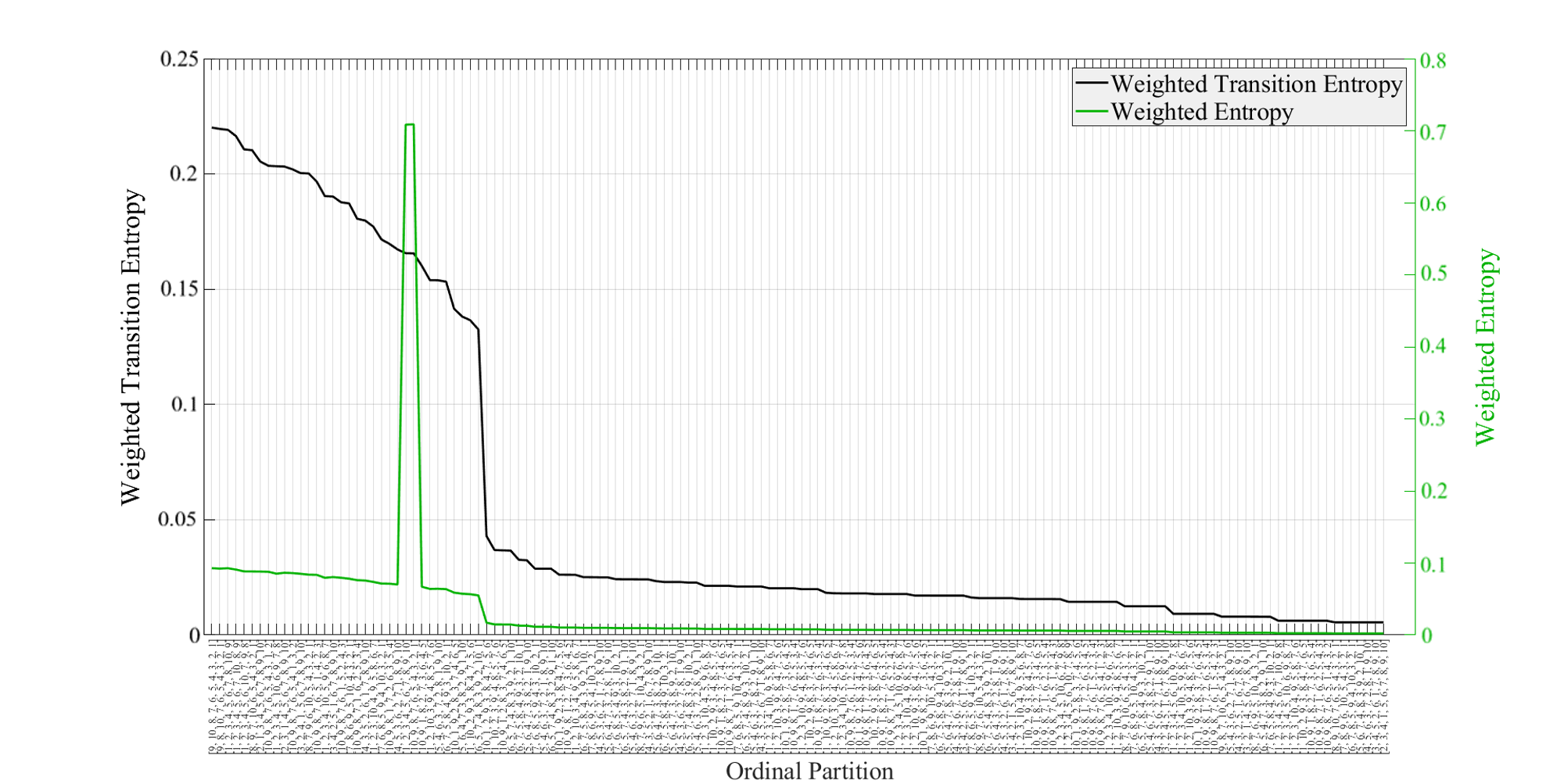}} \\ (a) \\
        {\includegraphics[width = 0.45\linewidth]{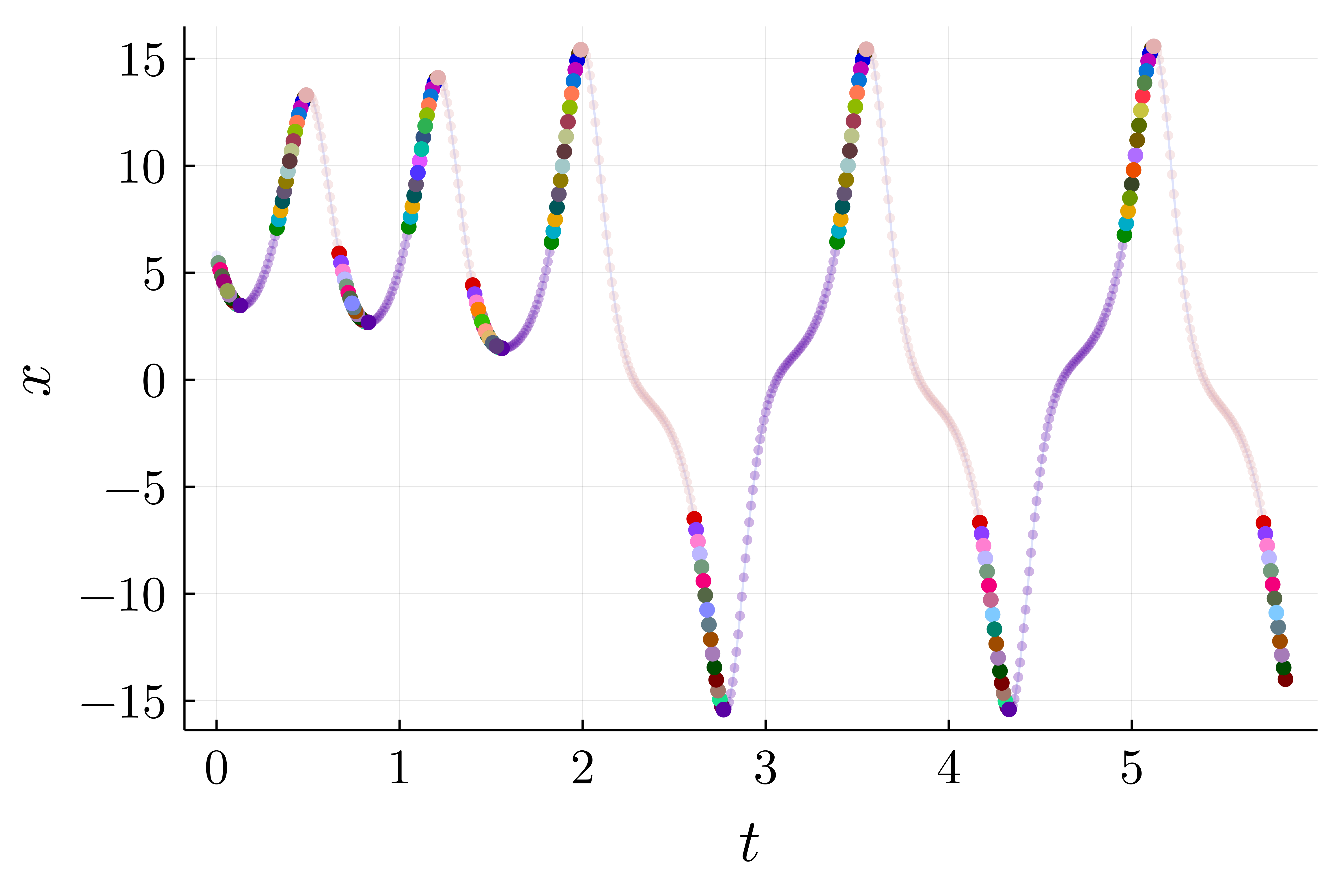}} \\ (b)
\end{tabular} 
&
\begin{tabular}{c} 
        {\includegraphics[width = 0.41\linewidth]{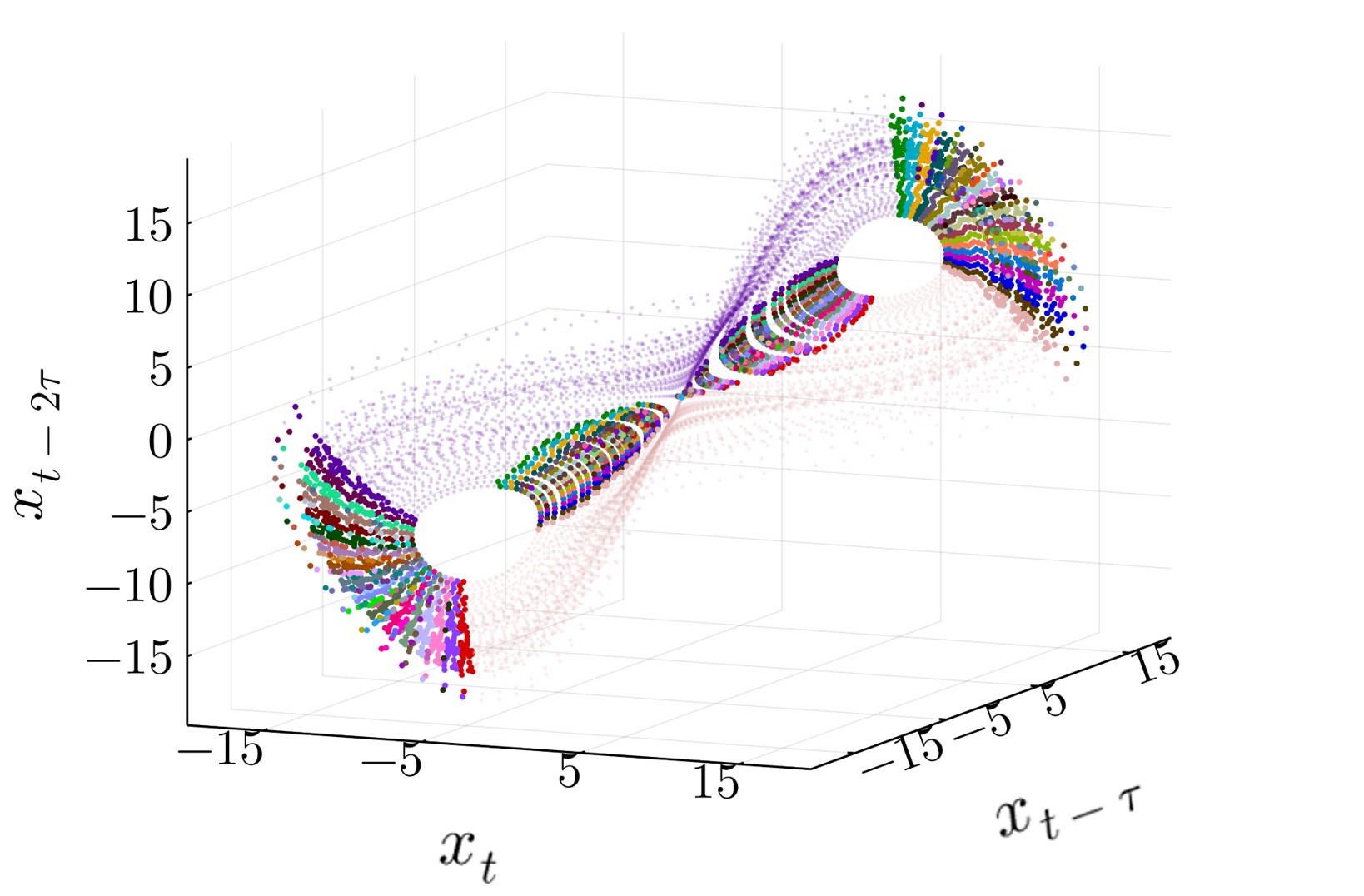}}\\ (c)\\
        {\includegraphics[width = 0.41\linewidth]{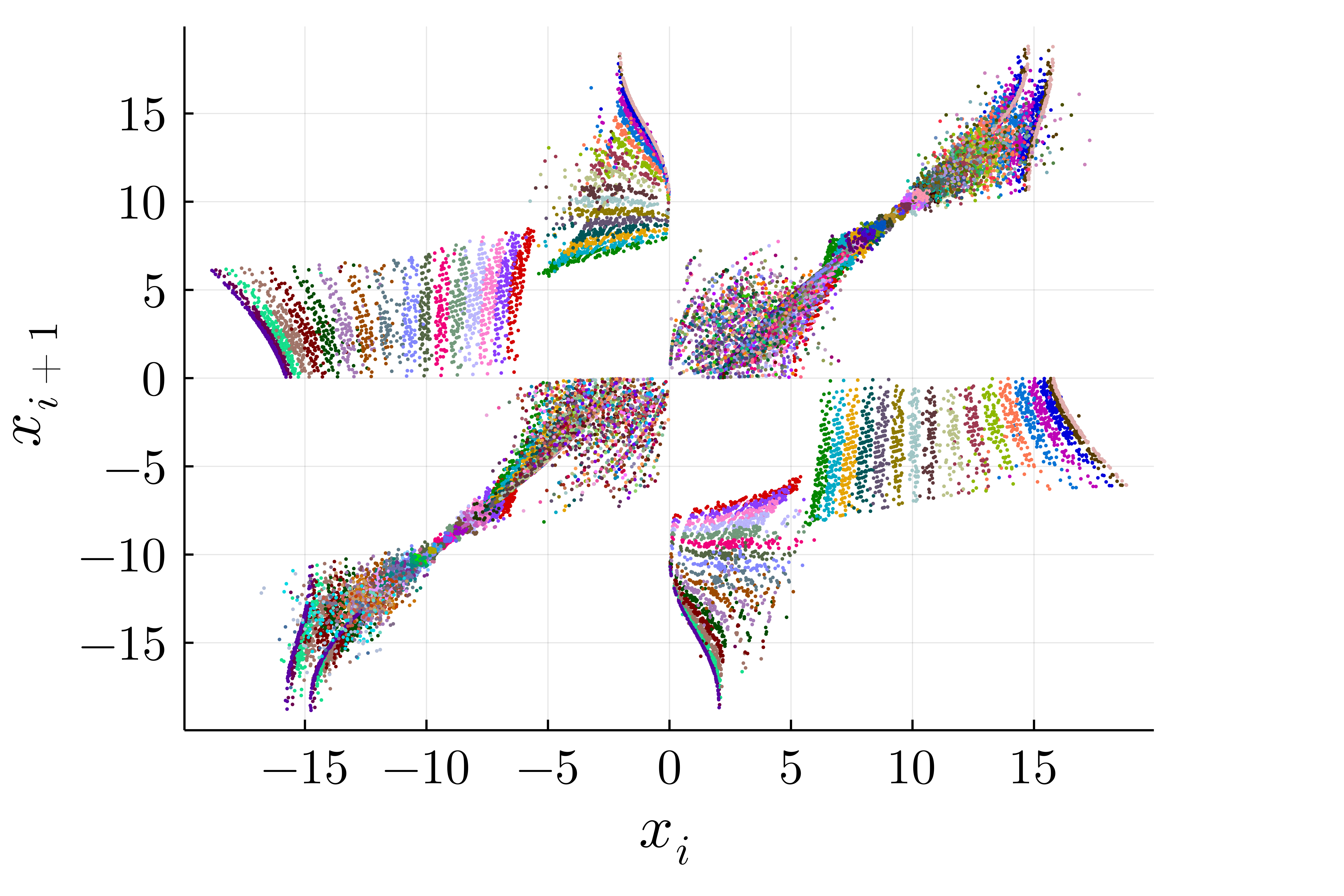}}\\ (d) 
\end{tabular}
\end{tabular}
    \caption{(a) The weighted transition entropy of each ordinal partition sorted in descending order (black line) and their corresponding weighted entropy (green line). (b) Lorenz's $x$ time series colored based on the ordinal partition. Each window's ordinal partition has been assigned to the first point of the window. (c) The Lorenz embedded attractor from its $x$-colored time series using Takens' embedding theorem with an embedding dimension of $M=3$ and an embedding lag of $T = 9$. The high entropy sections serve as appropriate Poincar\'e sections on the attractor. (d) Corresponding FRMs of each ordinal partition. N.B. Figure.~\ref{LorenzColoredOP_m10}c is not required to draw Fig.~\ref{LorenzColoredOP_m10}d, and it is only provided to demonstrate why high-entropy ordinal partitions can be considered good Poincar\'e sections.}
    \label{LorenzColoredOP_m10}
\end{figure*}

Figure~\ref{LorenzColoredOP_m10} depicts the results, equivalent to Fig.~\ref{SortedEntropy} and Fig.~\ref{LorenzColoredTransparentOP}, for the case when $m = 10$.

\subsection{\label{app:Rossler}R{\"o}ssler}

\begin{figure*}[htbp]
\begin{tabular}{cc}
\begin{tabular}{c}
        {\includegraphics[width = 0.4\linewidth]{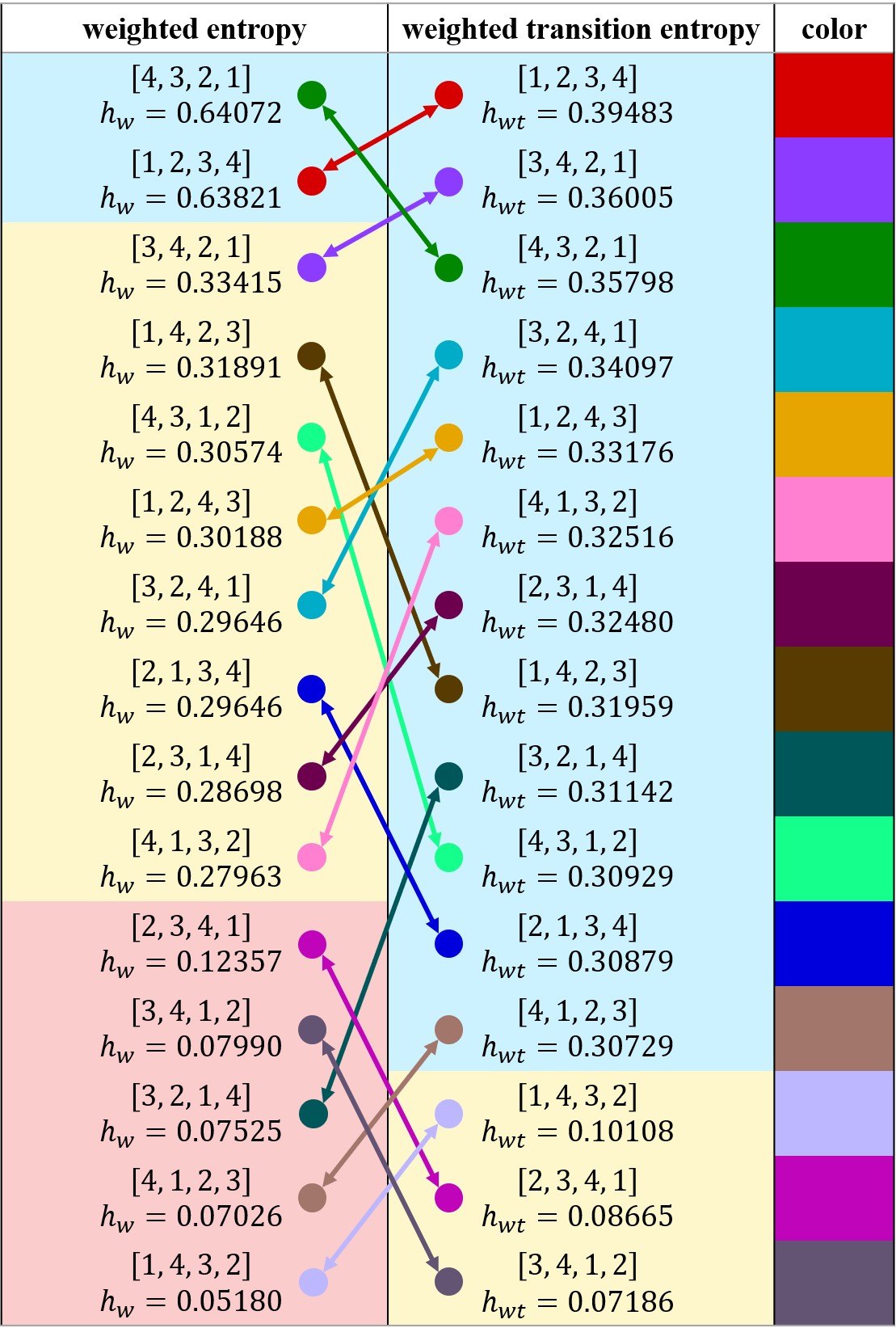}} \\ (a) \\
        {\includegraphics[width = 0.45\linewidth]{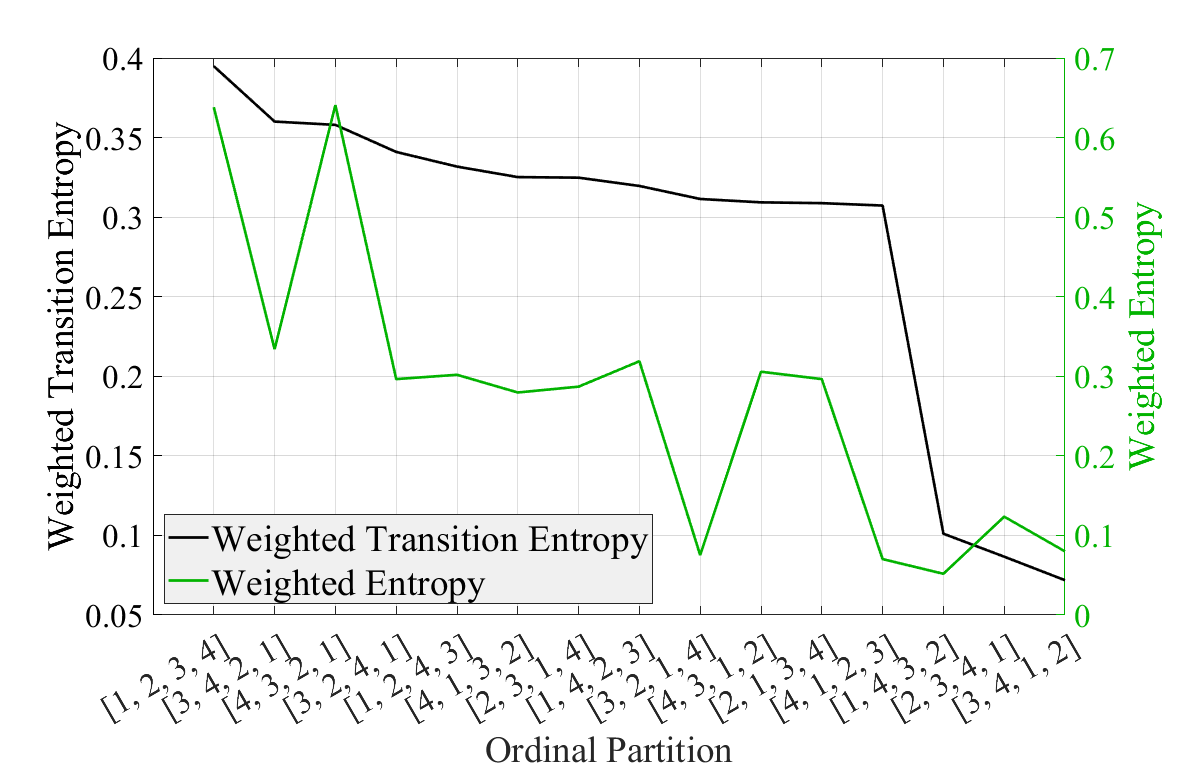}} \\ (b)
\end{tabular} 
&
\begin{tabular}{c} 
        {\includegraphics[width = 0.42\linewidth]{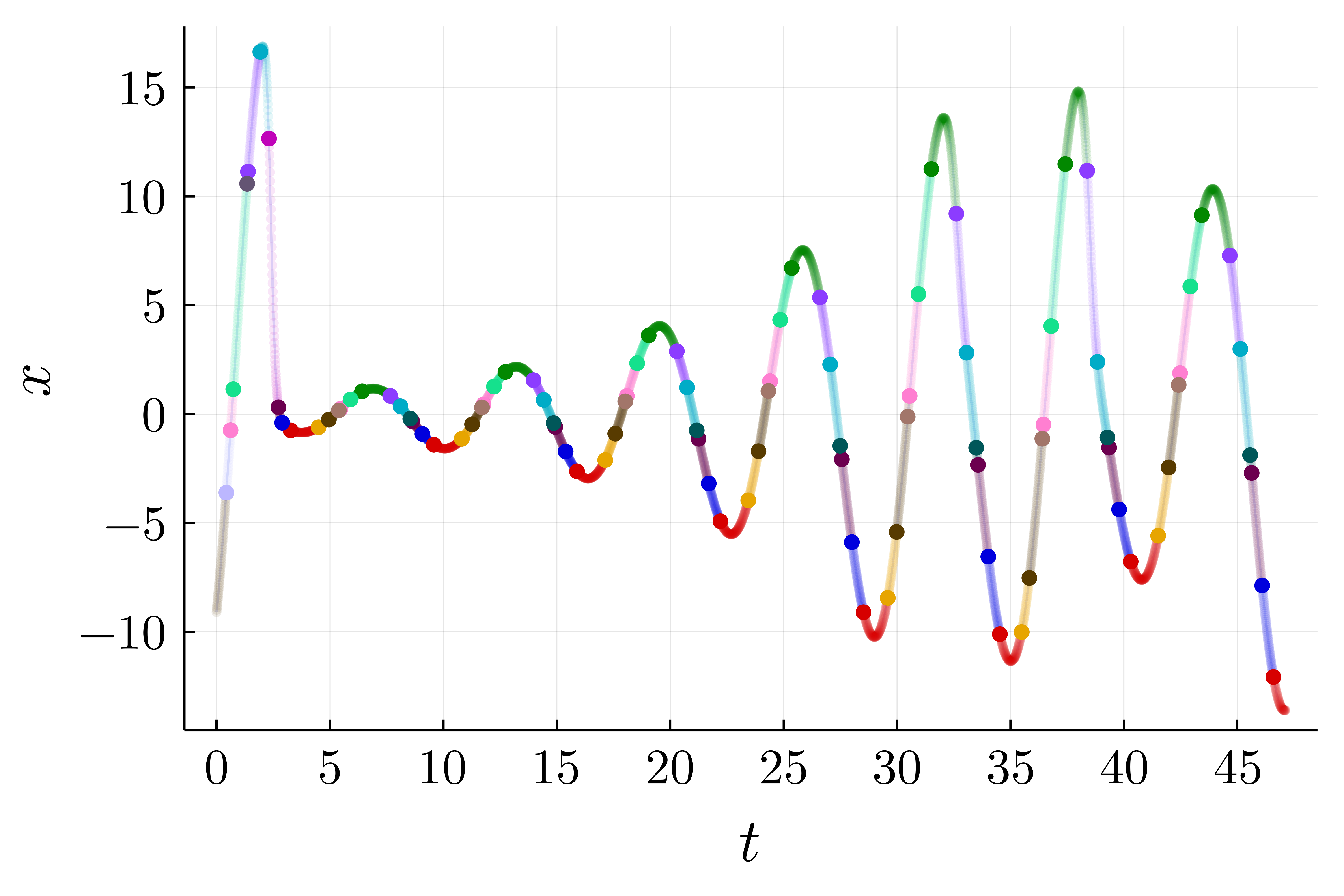}}\\ (c)\\
        {\includegraphics[width = 0.45\linewidth]{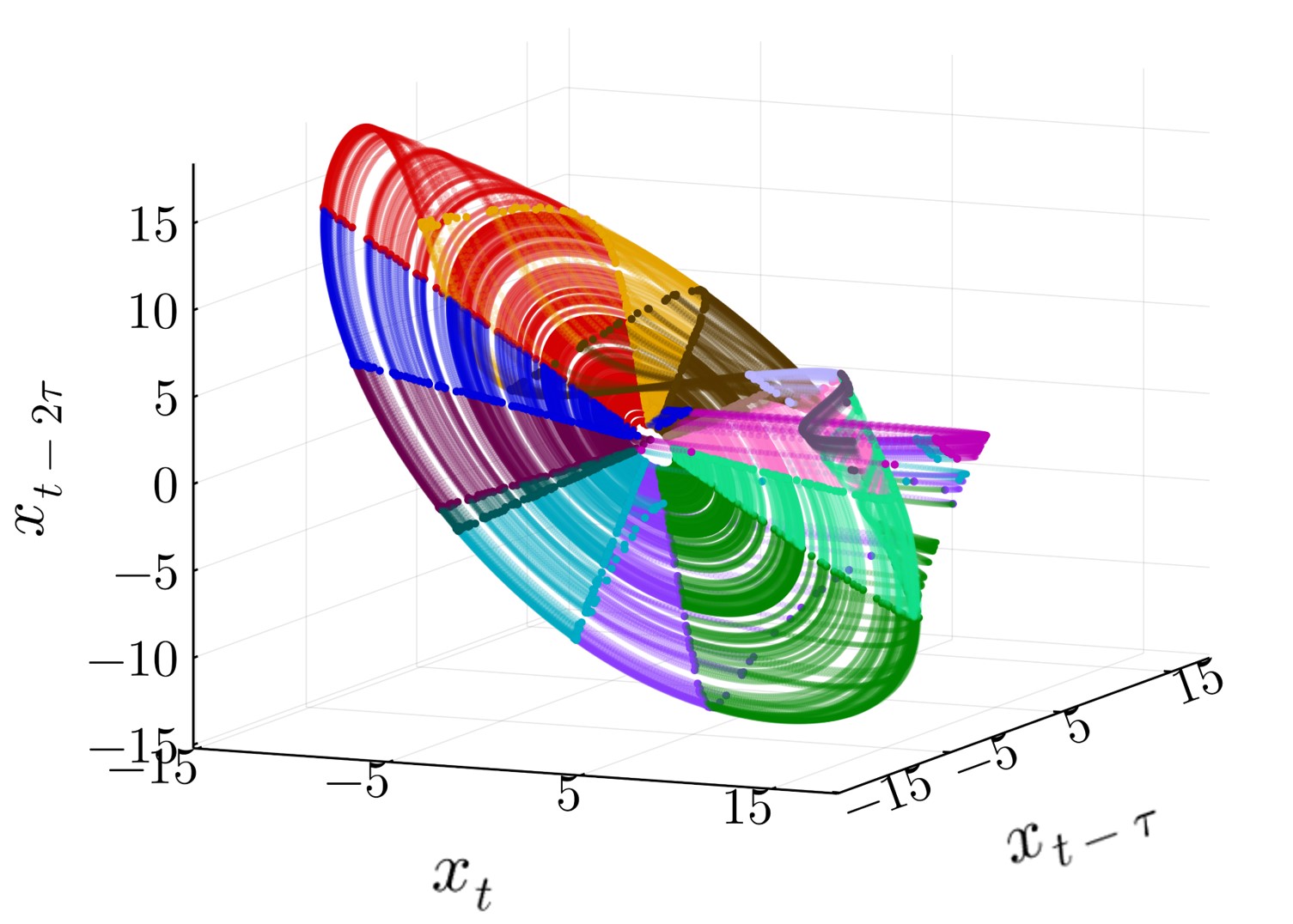}}\\ (d) \\
        {\includegraphics[width = 0.45\linewidth]{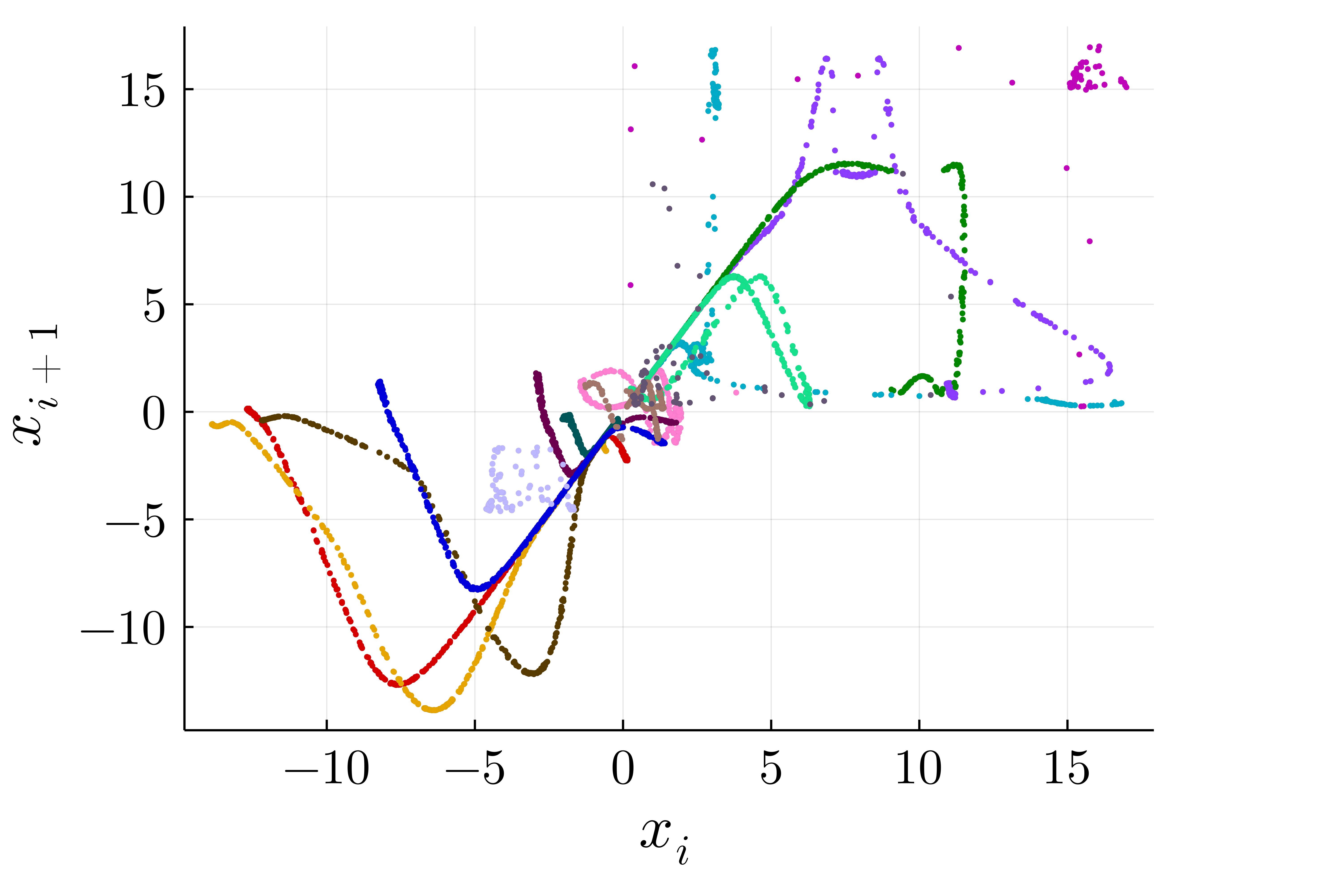}} \\ (e)
\end{tabular}
\end{tabular}
    \caption{(a) The list of ordinal partitions that have occurred on the R{\"o}ssler attractor. The first column is ordered by weighted entropy ($h_w$), and the second column is ordered by weighted transition entropy ($h_{wt}$). Each ordinal partition has been assigned a specific color, as shown in the third column. The background colors in the first and second columns represent the entropy level in their respective groups (See Fig.~\ref{RosslerColoredOP}b). Blue cells are in the first level of entropy, yellow cells are in the second level, and red cells are in the third level. There are three levels of weighted entropy ($h_w$) but only two levels of weighted transition entropy ($h_{wt}$). (b) The weighted transition entropy of each ordinal partition sorted in descending order (black line) and their corresponding weighted entropy (green line). (c) R{\"o}ssler's $x$ time series colored based on the ordinal partition. Each window's ordinal partition has been assigned to the first point of the window. (d) The R{\"o}ssler embedded attractor from its $x$-colored time series using Takens' embedding theorem with an embedding dimension of $M=3$ and an embedding lag of $T = 144$. The high entropy sections serve as appropriate Poincar\'e sections on the attractor. (e) Corresponding FRMs of each ordinal partition. N.B. Figure.~\ref{RosslerColoredOP}d is not required to draw Fig.~\ref{RosslerColoredOP}e, and it is only provided to demonstrate why high-entropy ordinal partitions can be considered good Poincar\'e sections.}
    \label{RosslerColoredOP}
\end{figure*}

Figure~\ref{RosslerColoredOP} illustrates the results of the application of our method --- equivalent to the results presented in Fig.~\ref{LorenzColoredOP}, Fig.~\ref{SortedEntropy} and Fig.~\ref{LorenzColoredTransparentOP} --- for the R{\"o}ssler system.
The equations of the R{\"o}ssler are as follows:

\begin{equation}
\label{Rossler}
\begin{array}{lll}
\cfrac{dx}{dt} & = & -y -z, \\
\cfrac{dy}{dt} & = & x + \alpha y, \\
\cfrac{dz}{dt} & = & \beta + (x - \gamma) z
\end{array}
\end{equation}
where $\alpha = 0.2$, $\beta = 0.2$, and $\gamma = 9$.

\subsection{\label{app:Mackey–Glass}Mackey–Glass}

\begin{figure*}[htbp]
\begin{tabular}{cc}
\begin{tabular}{c}
        {\includegraphics[width = 0.3\linewidth]{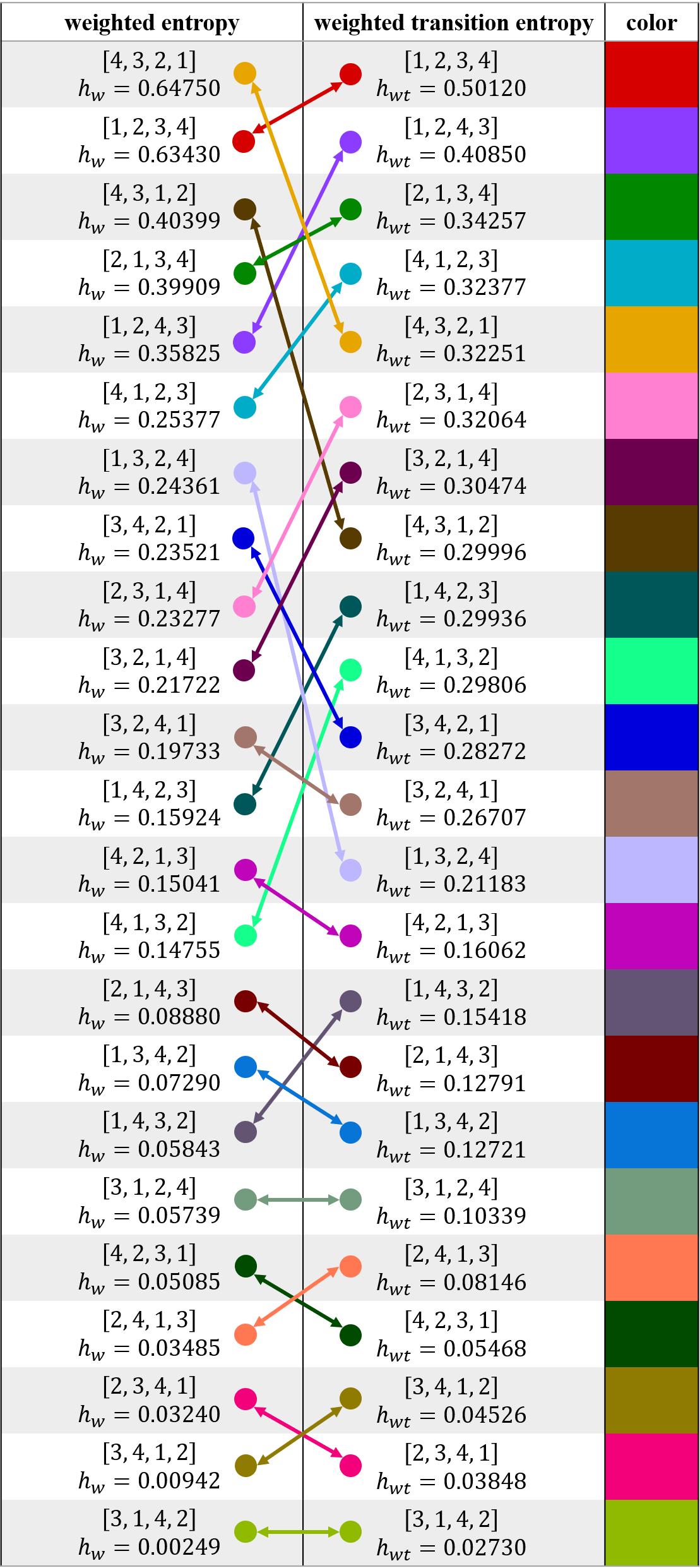}} \\ (a) \\
        {\includegraphics[width = 0.45\linewidth]{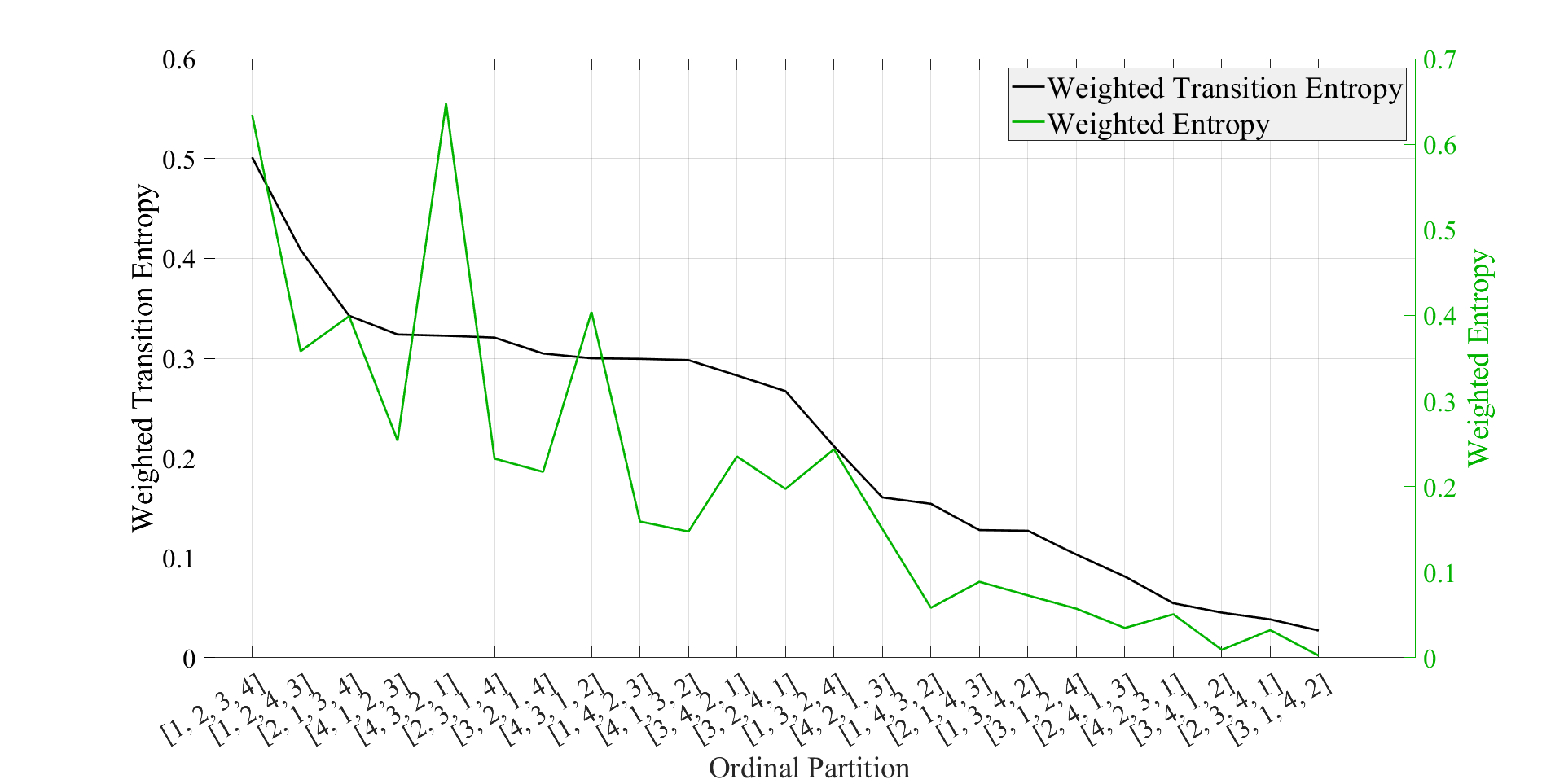}} \\ (b)
\end{tabular} 
&
\begin{tabular}{c} 
        {\includegraphics[width = 0.42\linewidth]{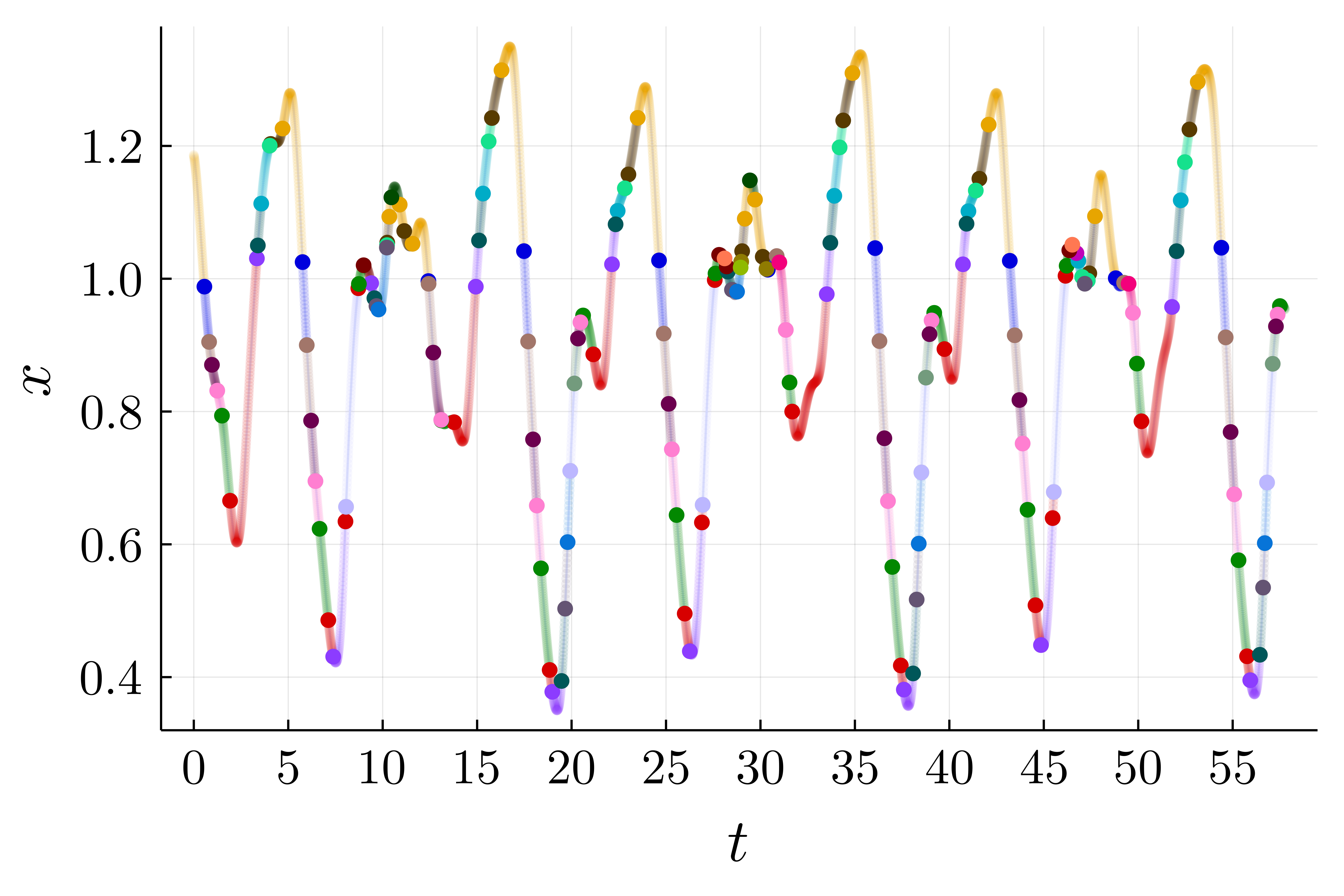}}\\ (c)\\
        {\includegraphics[width = 0.45\linewidth]{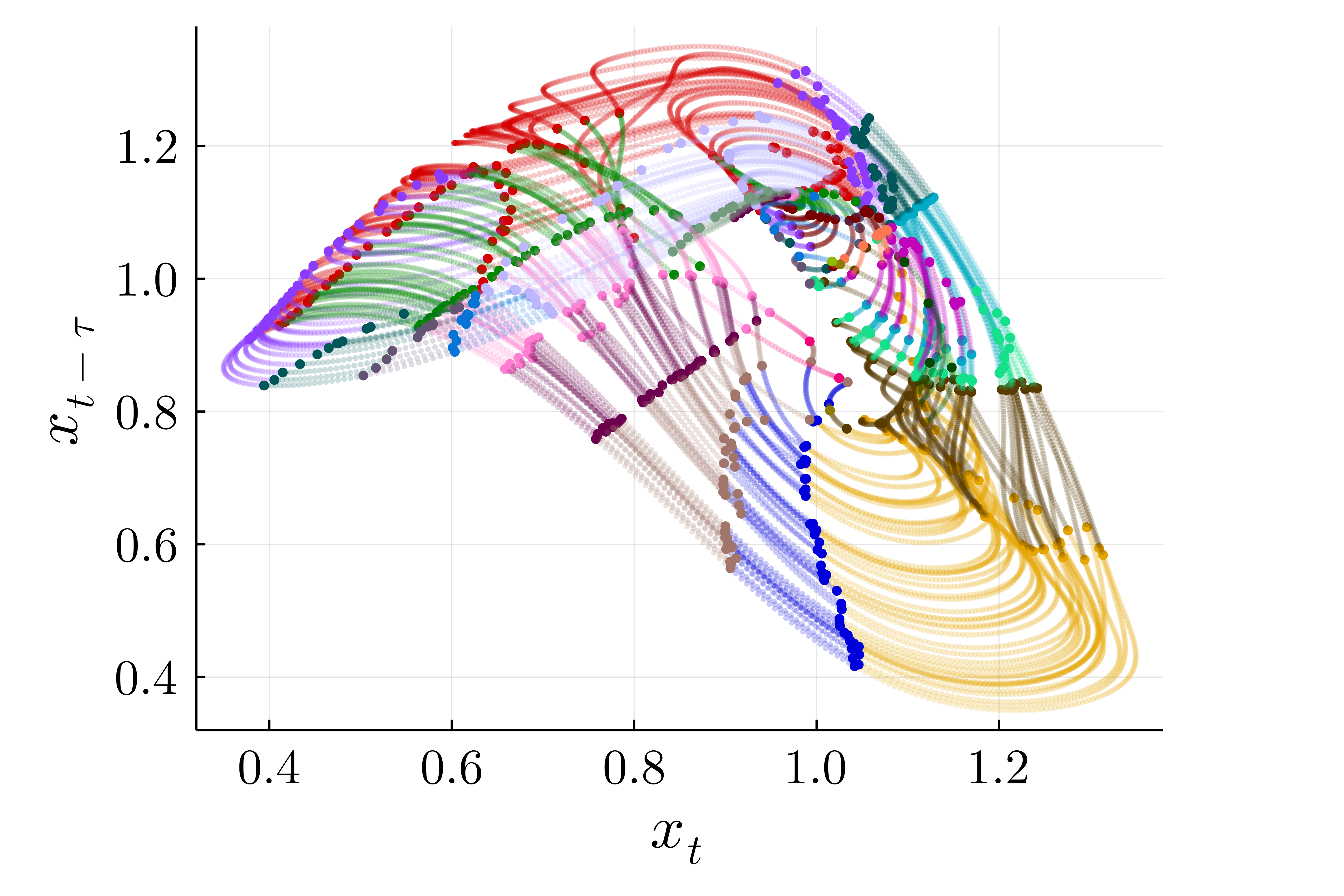}}\\ (d) \\
        {\includegraphics[width = 0.45\linewidth]{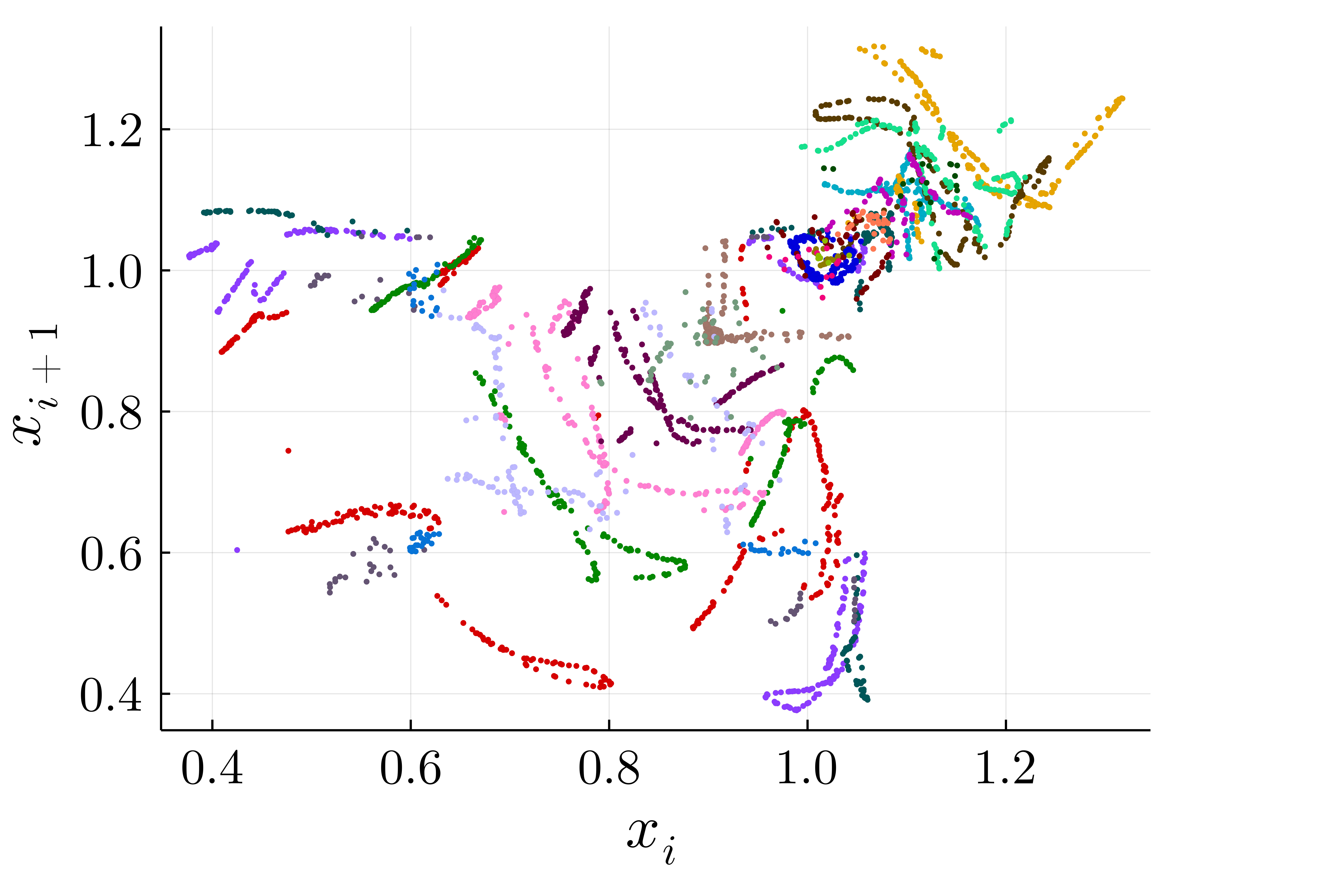}} \\ (e)
\end{tabular}
\end{tabular}
       \caption{(a) The list of ordinal partitions that have occurred on the MackeyGlass attractor. The first column is ordered by weighted entropy ($h_w$), and the second column is ordered by weighted transition entropy ($h_{wt}$). Each ordinal partition has been assigned a specific color, as shown in the third column. Because there are no obvious levels for different amounts of entropy (See Fig.~\ref{MackeyGlassColoredOP}b), the background color was not set for the first and second columns. (b) The weighted transition entropy of each ordinal partition sorted in descending order (black line) and their corresponding weighted entropy (green line). (c) MackeyGlass $x$ time series colored based on the ordinal partition. Each window's ordinal partition has been assigned to the first point of the window. (d) The MackeyGlass delayed embedded attractor with an embedding dimension of $M=2$ and an embedding lag of $T = 204$. The high entropy sections serve as appropriate Poincar\'e sections on the attractor. (e) Corresponding FRMs of each ordinal partition. N.B. Figure.~\ref{MackeyGlassColoredOP}d is not required to draw Fig.~\ref{MackeyGlassColoredOP}e, and it is only provided to demonstrate why high-entropy ordinal partitions can be considered good Poincar\'e sections.}
    \label{MackeyGlassColoredOP}
\end{figure*}

Figure~\ref{MackeyGlassColoredOP} shows the results equivalent to Fig.~\ref{LorenzColoredOP}, Fig.~\ref{SortedEntropy} and Fig.~\ref{LorenzColoredTransparentOP} for the Mackey-Glass system.
The equation of the Mackey-Glass is as follows:

\begin{equation}
\label{Mackey–Glass}
\cfrac{dx}{dt} = \beta \cfrac{x_{\tau}}{1 + {x_{\tau}}^n} - \gamma x 
\end{equation}
where $\gamma = 1$, $\beta = 2$, $\tau = 2$, $n = 9.65$, and $x_{\tau}$ represents the value of the variable $x$ at time $(t-\tau)$. 

\nocite{*}
\bibliography{sample}

\end{document}